\documentclass[10pt,fullpage]{article}
\usepackage{enumerate}
\usepackage{pgf}
\usepackage{tikz}
\usetikzlibrary{arrows}
\usepackage[utf8]{inputenc}

\newtheorem{theorem}{Theorem}
\newtheorem{lemma}{Lemma}
\newtheorem{corollary}{Corollary}
\newtheorem{proposition}{Proposition}

\begin{document}
\title{Indecomposables live in all smaller lengths}

\author{Klaus Bongartz\\Universität Wuppertal\\Germany\thanks{E-mail:bongartz@math.uni-wuppertal.de}}

\date{ Dedicated to A.V.Roiter and P.Gabriel}
\maketitle
\begin{abstract}
 We show that there are no gaps in the lengths of the indecomposable modules of a finite dimensional algebra over an algebraically closed field.
This result extends to indecomposables in any k-linear abelian category where the endomorphism algebras of the simples are k. For the proof we show that any distributive minimal representation-infinite algebra is isomorphic to its ray category and it has an interval-finite universal cover with a free fundamental group.
\end{abstract}

\section*{Introduction}

The  meaning of the title is the following result:

\begin{theorem} Let $A$ be an associative algebra of finite dimension over an algebraically closed field k. If there is a non-simple indecomposable $A$-module of length $n$, there is also one of length $n-1$.
\end{theorem} 

For representation-finite algebras - i.e. algebras having only finitely many isomorphism classes of finite-dimensional indecomposable representations - one even knows  since a long time that any non-simple module is an extension of a simple and an indecomposable. This is probably not true in general although Ringel has recently shown that any non-simple indecomposable is an extension of indecomposables.

The statement of the theorem is very naive and elementary, but the proof given here is not. It depends strongly on the work of Roiter, Gabriel and others on algebras of finite representation type. This is summarized very well in  chapters 13 and 14 of their book about representations of finite dimensional algebras. To describe in this introduction shortly the simple strategy of the proof I freely use notions from that book without any further explanation. Later on, there will be precise references.

 Up to Morita-equivalence one can assume that all simple modules have dimension one thereby replacing the length by the dimension.  In view of the known representation-finite
case it only remains to be shown that a minimal representation-infinite algebra - i.e. the algebra itself is not representation-finite, but each proper quotient is - admits indecomposables in all dimensions. If the algebra is distributive and zigzag-free we appeal to the well-elaborated covering theory and we can complete the proof by 'standard' arguments. So all we still have to deal with are non-distributive algebras or distributive algebras with a zigzag.

For non-distributive algebras there is a simple direct construction that was later on modified by Ringel to prove the existence of a so called accessible module in each dimension ( a first version of this article in the archive dates back to april 2009 ).

For distributive algebras with a zigzag I could not find a direct proof. The reason seems to be that the technique of cleaving diagrams due to Bautista, Larri\'{o}n and Salmer\'{o}n allows to construct easily infinite families of indecomposables, but their dimensions cannot be determined precisely. This situation is familiar from  modular representation theory: Induction implies  that the whole group has the same representation type as a p-Sylow subgroup, but the decomposition of an induced module into indecomposables can be very complicated. Thus for group algebras the two well-known Brauer-Thrall conjectures are easy to prove  whereas theorem 1 is not. This might be the reason why it was not formulated as the 'Brauer-Thrall 0 conjecture'. 

Now  to handle the distributive algebras with a zigzag I had  to generalize  a little bit a central result about coverings which might be  of independent interest:

\begin{theorem} Let $A$ be a distributive basic associative algebra with associated ray category $\vec{A}$. Suppose that $A$ is minimal  representation-infinite. Then we have:
\begin{enumerate}
 \item $\vec{A}$ has an interval-finite universal cover.
\item The fundamental group is free.
\item $A$ is isomorphic to the linearization $k\vec{A}$ of $\vec{A}$.
\end{enumerate}
\end{theorem} 

As I learned from Vossieck the last two statements of the theorem have already been obtained by Geiss in his unpublished diploma-thesis from 1990  by observing that Fischbachers arguments can be adapted to include the minimal representation-infinite case. However this does  not work for  the first statement and our inductive proof based on Fischbachers result gives all three statements at once.\vspace{0.8cm}

Theorem 1 has some nice corollaries. The first one is a generalization  from $A$-modules to objects of finite length in an abelian $k$-linear category. The anologous generalization is not true for the Brauer-Thrall conjectures as trivial examples show.
\begin{corollary}
Let $C$ be an abelian $k$-linear category over an algebraically closed field $k$. Suppose that all simple objects in $C$ have endomorphism algebra $k$. If there is an indecomposable non-simple object in $C$ of length $n$, there is also one of length  $n-1$.

\end{corollary}

Note that by a well-known counting argument the assumption on the endomorphism algebras is always true for modules over algebras whose dimension is strictly smaller than the cardinality of the field. So for example, the corollary applies to complex representations of Kac-Moody-algebras and their deformations.

\begin{corollary}( The naive criterion for finite representation type )
 The following conditions are equivalent for an algebra $A$ of finite dimension over an algebraically closed field.
\begin{enumerate}
 \item $A$ is representation-finite.
 \item There is a natural number $n$ such that there is no indecomposable $A$-module of that length.
\end{enumerate}
In fact, under these conditions the number $n= 2 \cdot dim A + 1000$ will always do.
\end{corollary}
 
Back in 1974 I tried to finish my diploma-thesis by applying this criterion. Much to my surprise Gabriel rejected my 'solution' because that obvious criterion was not proven. Now it is - hopefully.

The article is organized as follows. In chapter 1 we consider non-distributive algebras and in the central chapter 2 we study crowns - i.e. periodic zigzags - in minimal representation-infinite ray categories where at least one composition of irreducible morphisms does not vanish. The proof of the main reduction  resembles the proofs of some of the main results in the article on multiplicative bases. As explained at the beginning of section 2.2 it requires  first a finite strategy and second enough energy to carry this through. Finally in chapter 3 the statements made in this introduction are easy consequences of the general theory developed for representation-finite algebras.

I want to thank Dieter Vossieck for drawing my attention to an error in lemma 9 as stated in the first archive version.

For the sake of simplicity we will always work over an algebraically closed field $k$ of arbitrary characteristic, but the only thing that we really need is that  all simple modules have endomorphism algebra $k$. We consider left modules.\newpage

\section{The non-distributive case}

Now $A$ denotes  a basic associative algebra  of finite dimension over $k$ with Jacobson radical $J$. Such an algebra is given by a uniquely determined quiver $Q$ and a two-sided ideal $I$ inside the path algebra $kQ$, that is generated  by certain linear combinations of paths of length $\geq 2$. There is a commutative semi-simple subalgebra $B$ in $A$ that is a vector space supplement of $J$. For the next proof we need the following easy observation.

\begin{lemma}\label{konstrukt}
 Let $M$ be a non-zero $A$-module with an endomorphism $\phi$. Then we have:
\begin{enumerate}
 \item $\phi$ is nilpotent iff the induced map on the socle of $M$ is nilpotent.
 \item $\phi$ is nilpotent iff the induced map on the top $M/JM$ is nilpotent.
\item $M$ is indecomposable iff each endomorphism has exactly one eigenvalue.
\end{enumerate}

\end{lemma}

Recall that $A$ is distributive if its ideal lattice is distributive. This is equivalent to the fact that for all primitive idempotents $e,f$ the algebra $eAe$ is uniserial and that $fAe$ is cyclic as an $fAf$ left module or as an $eAe$ right module ( \cite[13.2]{Buch} ). Thus, if $A$ is not distributive, there are not necessarily different primitive idempotents $e,f$ and a natural number $l$ such that for the radical filtration $(R^{i})$ of $fAe$ as a $fAf-eAe$-bimodule we have $dim  R^{i}/R^{i+1} =1$ for all $i<l$, but $dim  R^{l}/R^{l+1} \geq 2.$ We choose elements $v, w$ in $R^l$ whose images in  $R^{l}/R^{l+1}$ are linearly independent, and we look at the two-sided ideal $K$ of $A$ generated by $R^{l+1}, Jv,vJ,Jw,wJ$. In the quotient $A/K$ we obtain primitive idempotent elements $\overline{e},\overline{f}$ and linearly independent elements $\overline{v},\overline{w}$ that are annihilated on both sides by the Jacobson radical of the quotient. Our aim is to construct in each dimension an indecomposable $A/K$-module. To simplify the notation a little bit we assume right from the beginning that the original $e,f,v,w$ have the properties mentioned before.

Let $d$ be the dimension of the indecomposable projective $Ae$. All non-zero quotients of this local module are again local, whence indecomposable. Thus we easily find indecomposables of dimension $m$ for all $m\leq d$. In particular, the family $Ae/\langle v-xw \rangle,$ $x \epsilon k,$ consists of pairwise non-isomorphic indecomposables. Here and later on we denote by $\langle X \rangle $ the $k$-subspace generated by some subset $X$ inside some vector space. 

 To construct at least one indecomposable in each dimension  we take the Kronecker-modules as our proto-types.
So let $n\geq 2$ be a natural number and take $n$ copies $Ax_{1},Ax_{2},\ldots ,Ax_{n}$ of the indecomposable projective $Ae$. Define $M=\oplus _{i=1}^{n} Ax_{i}$ and introduce the two subspaces $U_{0}=\langle wx_{i}- vx_{i+1}\mid 1\leq i \leq n-1 \rangle$ and $H= \langle vx_{1},vx_{2},\ldots ,vx_{n}\rangle$. Note that these are actually semi-simple submodules because $v$ and $w$ are annihilated by $J$.

\begin{lemma}\label{nondis} Using all the notations from above, the following is true:
\begin{enumerate}
 \item Let $U$ be a submodule of $M$  such that
\begin{enumerate}[i)]
  \item $U$ contains $U_{0}$ and also $wx_{n}$ for $U\neq U_{0}$.
  \item $U$ is contained in the radical $JM$ of $M$.
  \item $U \cap H =0$.
\end{enumerate} Then $M/U$ is indecomposable.
\item $U_{0}$ satisfies the conditions in part a). Let $U$ be a fixed maximal submodule satisfying these conditions. Then the socle of $N=M/U$ is isomorphic to $H$ under the canonical projection $\pi:M \longrightarrow M/U$.
\item Let $V$ be a  submodule of $N$ that contains the submodule $V_{0}$ generated by $\pi x_{1},\pi x_{2},\ldots ,\pi x_{n-1}$ and is contained in $V_{0}+ JN$. Then $V$ is indecomposable.
\item For each $m$ with $nd-(n-1) \geq m \geq (n-1)d-(n-2)$ we find an indecomposable subquotient of $M$ with dimension $m$.
\end{enumerate}

\end{lemma}

Proof: a) Condition i) just says $\overline{vx_{i+1}}=\overline{wx_{i}}$ for $1\leq i \leq n-1$ in $\overline{M}=M/U$. Here $\overline{x}$ denotes as usual  the image of an element of $M$ in $M/U$. Condition ii) guarantees that the top of $M/U$ is  $\oplus_{i=1}^{n} k\overline{x_{i}}$. Finally, the elements $\overline{vx_{1}},\overline{vx_{2}},\ldots ,\overline{vx_{n}}$ are still linearly independent in $\overline{M}$ by iii).

We claim that for all $i$ with $0\leq i \leq n-1$ we have $$(v^{-1}w)^{i}(\overline{M})= \langle \overline{x_{i+1}},\overline{x_{i+2}}, \ldots ,\overline{x_{n}}\rangle + J\overline{M}.$$
The start of the induction is trivial. In the step from $(i-1)$ to $i$ we get $w(v^{-1}w)^{i-1}(\overline{M})=\langle \overline{wx_{i}}, \ldots ,\overline{wx_{n}}\rangle$ because of $wJ=0$. Here $\overline{wx_{n}}=0$ for $U\neq U_{0}$. Taking the inverse image under multiplication with $v$ we find after  a short calculation $$(v^{-1}w)^{i}(\overline{M})= v^{-1}\langle \overline{wx_{i}}, \ldots ,\overline{wx_{n}}\rangle=\langle \overline{x_{i+1}},\overline{x_{i+2}}, \ldots ,\overline{x_{n}}\rangle + J\overline{M}.$$

Take now any endomorphism $\phi$ of $\overline{M}$. Then $\phi$ respects the flag 
$$\overline{M} \supseteq (v^{-1}w)\overline{M}\supseteq \ldots \supseteq (v^{-1}w)^{n-1}\overline{M} \supseteq 0.$$
Therefore we have $$\phi(\overline{x_{i}})= \sum_{j=i}^{n} \phi_{ij}\overline{x_{j}} +r_{i}$$ for some appropriate scalars $\phi_{ij}$ and $r_{i}$ in $J\overline{M}$. We obtain $$\phi(\overline{wx_{i}})= \sum_{j=i}^{n} \phi_{ij}\overline{wx_{j}}=\phi(\overline{vx_{i+1}})= \sum_{j=i+1}^{n} \phi_{i+1j}\overline{vx_{j}}$$ for all $1\leq i \leq n-1$,
 whence $ \phi_{ii}\overline{wx_{i}}= \phi_{ii}\overline{vx_{i+1}}=\phi_{i+1i+1}\overline{vx_{i+1}}$. Thus we have $\phi_{ii}=\phi_{11}=:a$ for all $i$. Therefore, $\phi - a\cdot id$ induces a nilpotent endomorphism on the top of $\overline{M}$ and also on $\overline{M}$ by the last lemma. So $a$ is the only eigenvalue of $\phi$ and  $\overline{M}$ is indecomposable.

b) It is easy to see that $U_{0}$ satisfies all conditions. Let $U$ be a maximal submodule with that property. Under the projection $\pi: M \longrightarrow M/U$ the semisimple module $H$ is embedded into the socle of $N$. We have to prove that any simple submodule of $N$ lies in the image of $H$. If $S$ is a simple submodule of $N$ its inverse image $I=\pi^{-1}S$ contains $U$ properly so that $I$ cannot satisfy all three conditions of part a). 

The first condition holds for $I$. So assume $I$ is not contained in $JM$. Then $S$ is not contained in $JN$ because of $\pi^{-1}(JN)=JM$. Using $S\cap JN = 0$ we can choose a $B$-module supplement $N'$ of $S$ in $N$ that contains $JN$. Then we get $N=S\oplus N'$ even as an $A$-module. But $N$ is indecomposable by part a). We conclude $n=1$ contradicting our assumption $n\geq 2$.

 Therefore $I$ cannot satisfy the third condition, i.e. $T=I \cap H \neq 0$. This implies 
$$U \subset \pi^{-1}(S)=U+T \subseteq U+H$$ and therefore $S=\pi(\pi^{-1}(S)) \subseteq \pi(H)$.

c) Of course the socle of $V$ contains $\overline{vx_{1}},\ldots \overline{vx_{n-1}}$ and also $\overline{vx_{n}} =\overline{wx_{n-1}}$, i.e. $\pi(H)$ which in turn is the socle of $N$. Since $V$ is a submodule of $N$, the socle of $V$ is $\pi(H)$. Write $K$ for the supspace of $V$ consisting of all elements killed by multiplication with $v$ and $w$. So $K$ contains the radical $JV$, but it will be strictly bigger in general. We claim that  $$(v^{-1}w)^{i}(V)= \langle \overline{x_{i+1}},\overline{x_{i+2}}, \ldots ,\overline{x_{n-1}}\rangle + K$$ holds for all $i=0,1, \ldots n-2.$ This is true for $i=0$ because $V$ is contained in $V_{0}+ JN$.  The induction-step is easy and similar to that in part a).

Now take an endomorphism $\phi$ of $V$. It respects the filtration by the  $(v^{-1}w)^{i}(V)$ and we have  $$\phi(\overline{x_{i}})= \sum_{j=i}^{n-1} \phi_{ij}\overline{x_{j}} +r_{i}$$ for some appropriate scalars $\phi_{ij}$ and $r_{i}$ in $K$. We obtain $$\phi(\overline{wx_{i}})= \sum_{j=i}^{n-1} \phi_{ij}\overline{wx_{j}}=\phi(\overline{vx_{i+1}})= \sum_{j=i+1}^{n-1} \phi_{i+1j}\overline{vx_{j}}$$ for all $i$ with $1\leq i \leq n-2$,
 whence $ \phi_{ii}\overline{wx_{i}}= \phi_{ii}\overline{vx_{i+1}}=\phi_{i+1i+1}\overline{vx_{i+1}}$. Thus we have $\phi_{ii}=\phi_{11}=:a$ for all $i\leq n-1$ and in addition $\phi(\overline{vx_{n}})=\phi(\overline{wx_{n-1}})=a\overline{vx_{n}}$. Therefore, $\phi - a\cdot id$ induces a nilpotent endomorphism on the socle of $V$  and also on $V$ by the last lemma. So $a$ is the only eigenvalue of $\phi$ and $V$ is indecomposable.

d) Choose a complete flag of submodules between $U_{0}$ and the maximal module $U$ fixed in part d). Dividing $M$ by these modules produces indecomposables of dimensions $m$ with $dim M/U_{0}= nd - (n-1) \geq m \geq dim M/U$. Similarly a complete flag of submodules of $N=M/U$ starting with $V_{0}$ and ending with $V_{0}+JN$ gives us indecomposables with dimensions ranging from $dim V_{0}$ to $dim M/U -1$. Since $V_{0}$ is generated by $\overline{x_{1}}, \ldots \overline{x_{n-1}}$ we have $dim V_{0} \leq (n-1)d-(n-2)$.\vspace{0.8cm}

The reader is invited to verify what the preceding construction means at least in the following two simple, but typical examples that are given by quivers with relations.\vspace{0.3cm}

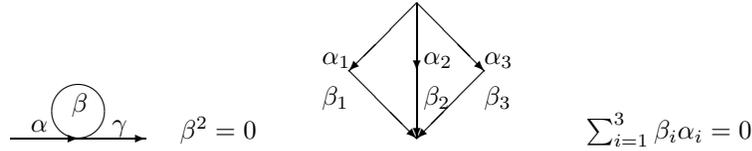
\setlength{\unitlength}{0.9cm}
\begin{picture}(10,2.5)
\put(0,0){\vector(1,0){1}}
\put(1,0){\vector(1,0){1}}
\put(1,0.4){\circle{0.8}}

\put(6,2){\vector(1,-1){1}}
\put(6,2){\vector(0,-1){1}}
\put(6,2){\vector(-1,-1){1}}
\put(5,1){\vector(1,-1){1}}
\put(6,1){\vector(0,-1){1}}
\put(7,1){\vector(-1,-1){1}}

\put(0.3,0.1){$\alpha$}
\put(0.9,0.4){$\beta$}
\put(1.5,0.1){$\gamma$}
\put(2.5,0){$\beta^{2}=0$}

\put(4.6,1.1){$\alpha_{1}$}
\put(6.1,1.1){$\alpha_{2}$}
\put(7,1.1){$\alpha_{3}$}
\put(4.6,0.5){$\beta_{1}$}
\put(6.1,0.5){$\beta_{2}$}
\put(7,0.5){$\beta_{3}$}
\put(8.5,0){$\sum_{i=1}^{3} \beta_{i}\alpha_{i}=0$}
\end{picture}\\

Our findings can be summarized in the following result:
\begin{proposition}\label{nondisinf}
Let $A$ be a basic non-distributive algebra. Then there is an indecomposable $W$ of countable dimension having for each natural number $m$ an indecomposable subquotient $V$ of dimension $m$.
 
\end{proposition}

Proof: We use the reductions and notations introduced before. Let $M'$ be the direct sum of an infinite sequence $Ax_{i}$ of copies of $Ae$ and let $U'$ be the submodule generated by all differences $ wx_{i}- vx_{i+1}$. Clearly, there is an endomorphism $T$ of $M'$ that maps $x_{i}$ to $x_{i+1}$ for all i. Since $U'$ is $T$-invariant, we obtain an induced endomorphism $T$ on the quotient $W:=M'/U'$. 

To see that $W$ is indecomposable, we determine its andomorphism algebra $B$. As before one sees that any endomorphism $\phi$ respects the infinite descending chain of submodules $W_{i}$ generated by $JW$ and the $x_{j}$ with $j\geq i$. This implies inductively that an endomorphism $\phi$ with $\phi(x_{1}) \epsilon JW$ maps $W$ into $JW$. Let $I$ be the set of all endomorphisms $\phi$ with $\phi(W) \subseteq JW$. Then $I$ is a nilpotent ideal and we claim that $B$ is the direct sum of $I$ and the subalgebra $k[T]$ generated by $T$ which is isomorphic to the polynomial algebra in one indeterminate.

Indeed, let $\psi$ be an endomorphism. Then we have $\psi(x_{1})= \sum_{i=1}^{n}\lambda_{i}x_{i}+ j$ with some scalars $\lambda_{i}$ and $j$ in $JW$. Thus $\phi:=\psi - \sum_{i=1}^{n}\lambda_{i}T^{i}$ maps $x_{1}$ into $JW$, whence it belongs to $I$ and $B$ is the sum of $I$ and $k[T]$ whose intersection is trivial. It follows easily that $0$ and $1$ are the only idempotents in $B$. Therefore $W$ is indecomposable.

Finally, for a given natural number $m$ there is an $n$ with $nd-(n-1) \geq m \geq (n-1)d-(n-2)$. Then the module $M/U_{0}$ considered in lemma \ref{nondis} is a submodule of $W$ and it contains an indecomposable subquotient of dimension $m$.

\section{On crowns in minimal representation-infinite ray categories}
\subsection{Reminder on ray categories and cleaving diagrams}
Unfortunately, we have to recall now a lot of definitions and results mainly from the book \cite{Buch}.

A locally bounded category $k$-category $A$ is a $k$-linear category where different objects are not isomorphic, where all endomorphism algebras $A(x,x)$ are local and where the direct sums $\bigoplus_{y\epsilon A}A(x,y)$ and $\bigoplus_{y\epsilon A}A(y,x)$ are of finite dimension for all $x\epsilon A$. A finite dimensional $A$-module $M$ is a covariant $k$-linear functor from $A$ to the category of $k$-vectorspaces such that the sum of the dimensions of all $M(x), x \epsilon A$, is finite. $A$ is locally representation finite, if for any object $x\epsilon A$ there are up to isomorphism only finitely many indecomposable modules $U$ with $U(x)\neq 0$.

A locally bounded $k$-category is distributive  if all endomorphism algebras $A(x,x)$ are uniserial and all homomorphism spaces $A(x,y)$ are cyclic as an $A(x,x)$ right module or as an $A(y,y)$ left module. The product $A(x,x)^{\ast}\times A(y,y)^{\ast}$ of the two automorphism groups acts on $A(x,y)$ and the orbit of a morphism is the corresponding ray. These rays are the morphisms of the ray category $\vec{A}$ attached to $A$. The properties of $\vec{A}$ are subsumed in the following axioms that define the abstract notion of a ray category $P$ ( \cite[section 13.4]{Buch} ):\begin{enumerate}
 \item The objects form a set and they are pairwise not isomorphic.
\item There is a family of zero-morphisms $0_{xy}:x \rightarrow y $, $x,y \epsilon P$, satisfying $\mu 0=0=0\nu$ whenever the composition is defined.
\item For each $x \epsilon P$, $P(x,y)=\{0\}$ and $P(y,x)=\{0\}$ for almost all $y \epsilon P$.
\item For each $x$ one has $P(x,x)=\{ id_{x},\sigma,\ldots ,\sigma^{n-1}\neq 0=\sigma^{n}\}$. Here $n$ depends on $x$.
\item For each $x,y$, the set $P(x,y)$ is cyclic under the action of $P(x,x)$ or of $P(y,y)$.
\item If $\kappa, \lambda,\mu,\nu$ are morphisms with $\lambda \mu\kappa= \lambda \nu \kappa\neq 0$ then $\mu= \nu$.
 \end{enumerate}
 Starting with such an abstract ray category $P$ one constructs in a natural way its linearization $k(P)$, which is a locally bounded distributive $k$-category having the original category $P$ as the associated ray category $\vec{k(P)}$ ( \cite[ 13.5]{Buch} ). In sharp contrast, a locally bounded distributive category $A$ is in general not isomorphic to $k(\vec{A})$. If it is, $A$ is called standard. We say that $P$ is  ( locally ) representation finite or minimal representation-infinite if $k(P)$ is so, and this is independent of the field by \cite[ 14.7]{Buch}.

To study a ray category $P$ and its universal cover $\tilde{P}$ we look at the quiver $Q_{P}$ of $P$ ( \cite[section 13.6]{Buch} ). Its points are the objects of $P$ and its arrows the irreducible morphisms in $P$, i.e. those non-zero morphisms that cannot be written as a product of two morphisms different from identities. Each non-zero non-invertible $\mu$ in $P$ is then  a product of irreducible morphisms, and the depth $d(\mu)$ of $\mu$ is the maximal number of factors occuring in these products. The non-zero morphisms in $P$ are partially ordered by defining $\mu \leq \nu$ iff $\nu = \alpha \mu \beta $ for some morphisms $\alpha$ and $\beta$. A morphism is long if it is maximal with respect to this order and not irreducible. For $x \epsilon P$ we also consider the finite partially ordered set $x/P$ of all non-zero morphisms with domain $x$. Here we define $\phi \preceq \psi$ iff $\psi= \chi \phi$. The dual order on the set $P/x$ of the non-zero morphisms with a fixed codomain is also denoted by $\preceq$. The path category $\mathcal{P}Q_{P}$ has the points of $Q_{P}$ as objects and the paths in $Q_{P}$ as non-zero morphisms, to which we add formal zero-morphisms. There is a canonical full functor $\vec{}:\mathcal{P}Q_{P} \longrightarrow P$ from the path category to $P$ which is the 'identity' on objects, arrows and zero-morphisms. Two paths in $Q_{P}$ are interlaced if they belong to the transitive closure of the relation $R$ given by $(v,w)\epsilon R$ iff $v=pv'q, w=pw'q$ and $\vec{v}'=\vec{w}'\neq 0$ where $p$ and $q$ are not both identities. A contour of $P$ is a pair $(v,w)$ of non-interlaced paths with $\mu=\vec{v}=\vec{w} \neq 0$ ( see \cite[section 13.6]{Buch} ). Then we  say that $\mu$ occurs in the contour $(v,w)$. Note that these contours are called essential contours in  \cite{BGRS,Fischbacher}. A decomposition $v=v_{r}v_{r-1}\ldots v_{1}$ of a path is non-trivial if all subpaths $v_{i}$ have length $1$ at least. Similarly, a factorization of a morphism is non-trivial if none of the factors is an identity.

A functor $F:D \longrightarrow P$ between ray categories is cleaving ( \cite[ 13.8]{Buch} ) iff it satisfies the following two conditions and their duals: a) $F\mu =0$ iff $\mu=0$;
b) If $\alpha \epsilon D(x,y)$ is irreducible and $F\mu:Fx \rightarrow Fz$ factors through $F\alpha$ then $\mu$ factors already through $\alpha$.
The key fact about cleaving functors is that $P$ is not ( locally ) representation finite if $D$ is not.

In this article $D$ will always be given by its quiver $Q_{D}$, that has no oriented cycles, and some relations. Two  paths between the same points give always the same morphism, and zero relations are written down explicitely. As in \cite[section 13]{Buch} the cleaving functor is then defined by drawing the quiver of $D$ with relations and by writing the morphism $F\alpha$ in $P$ close to each arrow $\alpha$. To avoid confusions by to many letters in our figures we include sometimes not all names of morphisms ( see  figure 2.1 ) or we only mention all morphisms occurring in a figure in the text.
For instance, let $D$ be the ray category with the natural numbers as objects and with arrows $2n \leftarrow 2n+1$ and $2n+1 \rightarrow 2n+2$ for all $n$. Then a cleaving functor from $D$ to $P$ is called a zigzag in \cite[section 13.9]{Buch} and $P$ is said to contain a zigzag. A functor from $D$ to $P$ is just an infinite sequence of morphisms $(\sigma_{1},\rho_{1},\sigma_{2},\rho_{2},\ldots   )$ in $P$ such that $\rho_{i}$ and $\sigma_{i}$ always have common domain and $\rho_{i}$ and $\sigma_{i+1}$ common codomain.
 The  functor is cleaving iff none of the equations $\sigma_{i}=\xi \rho_{i}$,$\xi \sigma_{i}= \rho_{i}$,$\sigma_{i+1}\xi =\rho_{i}$ or $\sigma_{i+1}= \rho_{i} \xi$ has a solution. The situation is usually illustrated by the following zigzag:\vspace{0.5cm}

\setlength{\unitlength}{0.8cm}
\begin{picture}(20,2)

\put(1,2){\vector(1,-1){1}}

\put(1,2){\vector(-1,-1){1}}\put(-0.1,1.5){$\sigma_{1}$}

\put(3,2){\vector(-1,-1){1}}\put(0.9,1.5){$\rho_{1}$}

\put(3,2){\vector(1,-1){1}}\put(1.9,1.5){$\sigma_{2}$}
\put(5,2){\vector(-1,-1){1}}\put(2.9,1.5){$\rho_{2}$}

\put(5,2){\vector(1,-1){1}}
\put(7,2){\vector(-1,-1){1}}

\put(7,2){\vector(1,-1){1}}\multiput(9,1.5)(0.5,0){10}{-}
\put(6,0){figure 2.1}
\end{picture}\vspace{0.5cm}

A crown in $P$ of length $2n$  is a zig-zag that becomes periodic after n steps, i.e. one has  $\sigma_{i}=\sigma_{n+i}$ and $\rho_{i}=\rho_{n+i}$. If $P$ is finite and contains a zigzag, it contains also a crown.   We denote such a crown by $(\sigma_{1},\rho_{1},\sigma_{2}, \ldots ,\rho_{n})$. By axiom e) of a ray category the length of a crown is at least $4$.

Later on we need the following representation infinite ray categories. The numbers refer to the list in \cite[section 10.7]{Buch}.\vspace{1cm}\setlength{\unitlength}{0.7cm}

\begin{picture}(18,3)
\put(0,0){\begin{picture}(7,3)
 \put(0,3){\line(1,0){1}}
\put(1,3){\line(1,0){1}}
\put(2,3){\vector(1,0){1}}
\put(3,3){\line(1,0){1}}
\put(3,3){\vector(0,-1){1}}
\put(4,3){\line(1,0){1}}
\put(2,3){\vector(0,-1){1}}
\put(2,2){\vector(1,0){1}}\put(4,2){$12$}

\put(0,3){\circle*{0.1}}
\put(1,3){\circle*{0.1}}
\put(2,3){\circle*{0.1}}
\put(3,3){\circle*{0.1}}
\put(5,3){\circle*{0.1}}
\put(4,3){\circle*{0.1}}
\put(2,2){\circle*{0.1}}
\put(3,2){\circle*{0.1}}

\end{picture}}

\put(6,0)
{\begin{picture}(7,3)
 \put(0,2){\line(1,0){1}}
\put(1,2){\line(1,0){1}}
\put(2,3){\vector(1,0){1}}
\put(3,3){\line(1,0){1}}
\put(3,3){\vector(0,-1){1}}
\put(4,3){\line(1,0){1}}
\put(2,3){\vector(0,-1){1}}
\put(2,2){\vector(1,0){1}}
\put(0,2){\circle*{0.1}}
\put(1,2){\circle*{0.1}}
\put(2,3){\circle*{0.1}}
\put(3,3){\circle*{0.1}}
\put(5,3){\circle*{0.1}}
\put(4,3){\circle*{0.1}}
\put(2,2){\circle*{0.1}}
\put(3,2){\circle*{0.1}}\put(4,2){$14$}

\end{picture}}

\put(12,0){
\begin{picture}(7,3)

\put(1,3){\line(1,0){1}}
\put(2,3){\vector(1,0){1}}
\put(3,3){\line(1,0){1}}
\put(3,3){\vector(0,-1){1}}
\put(1,2){\line(1,0){1}}
\put(2,3){\vector(0,-1){1}}
\put(2,2){\vector(1,0){1}}
\put(1,3){\circle*{0.1}}
\put(2,3){\circle*{0.1}}
\put(3,3){\circle*{0.1}}
\put(1,2){\circle*{0.1}}
\put(4,3){\circle*{0.1}}
\put(2,2){\circle*{0.1}}
\put(3,2){\circle*{0.1}}\put(4,2){$11$}
\end{picture}}
\end{picture}\vspace{0.5cm}

\setlength{\unitlength}{0.7cm}
\begin{picture}(18,3)
\put(-1,0){
\begin{picture}(7,3)
\put(1,3){\line(1,0){1}}
\put(2,3){\vector(1,0){1}}
\put(3,3){\line(1,0){1}}
\put(3,3){\vector(0,-1){1}}
\put(4,3){\line(1,0){1}}
\put(5,3){\line(1,0){1}}
\put(6,3){\line(1,0){1}}
\put(2,3){\vector(0,-1){1}}
\put(2,2){\vector(1,0){1}}
\put(2,2){\circle*{0.1}}
\put(1,3){\circle*{0.1}}
\put(2,3){\circle*{0.1}}
\put(3,3){\circle*{0.1}}
\put(5,3){\circle*{0.1}}
\put(6,3){\circle*{0.1}}
\put(7,3){\circle*{0.1}}
\put(4,3){\circle*{0.1}}
\put(2,2){\circle*{0.1}}
\put(3,2){\circle*{0.1}}\put(5,2){$20$}

\end{picture}}

\put(8,0){
\begin{picture}(7,3)

\put(6,3){\vector(-1,0){1}}
\put(0,3){\line(1,0){1}}
\put(1,3){\line(1,0){1}}
\put(2,3){\line(1,0){1}}
\put(3,3){\line(1,0){1}}
\put(5,3){\vector(-1,0){1}}
\put(4,3){\vector(0,-1){1}}
\put(6,3){\vector(0,-1){1}}
\put(6,2){\vector(-1,0){2}}
\put(4,2){\circle*{0.1}}
\put(2,3){\circle*{0.1}}
\put(3,3){\circle*{0.1}}
\put(5,3){\circle*{0.1}}
\put(6,2){\circle*{0.1}}
\put(4,3){\circle*{0.1}}
\put(6,3){\circle*{0.1}}
\put(0,3){\circle*{0.1}}
\put(1,3){\circle*{0.1}}\put(2,2){$44$}

\end{picture}}
\end{picture}\vspace{0.5cm}

\setlength{\unitlength}{0.7cm}
\begin{picture}(18,3)
\put(0,0){
\begin{picture}(7,3)

\put(4,3){\vector(-1,0){1}}
\put(2,3){\line(1,0){1}}
\put(4,3){\vector(1,0){1}}
\put(5,3){\line(1,0){1}}
\put(3,3){\vector(0,-1){1}}
\put(4,3){\vector(0,-1){1}}
\put(5,3){\vector(0,-1){1}}
\put(4,2){\vector(1,0){1}}
\put(4,2){\vector(-1,0){1}}
\put(4,2){\circle*{0.1}}
\put(2,3){\circle*{0.1}}
\put(3,3){\circle*{0.1}}
\put(5,3){\circle*{0.1}}
\put(5,2){\circle*{0.1}}
\put(4,3){\circle*{0.1}}
\put(6,3){\circle*{0.1}}
\put(3,2){\circle*{0.1}}\put(6,2){$93$}

\end{picture}}

\put(8,0){
\begin{picture}(7,3)
\put(1,3){\line(1,0){1}}
\put(4,3){\vector(-1,0){1}}
\put(2,3){\line(1,0){1}}
\put(4,3){\vector(1,0){1}}
\put(5,2){\line(1,0){1}}
\put(3,3){\vector(0,-1){1}}
\put(4,3){\vector(0,-1){1}}
\put(5,3){\vector(0,-1){1}}
\put(4,2){\vector(1,0){1}}
\put(4,2){\vector(-1,0){1}}
\put(4,2){\circle*{0.1}}
\put(1,3){\circle*{0.1}}
\put(2,3){\circle*{0.1}}
\put(3,3){\circle*{0.1}}
\put(5,3){\circle*{0.1}}
\put(5,2){\circle*{0.1}}
\put(6,2){\circle*{0.1}}
\put(4,3){\circle*{0.1}}
\put(3,2){\circle*{0.1}}\put(2,2){$96$}\put(-1,1){figure 2.2}

\end{picture}}
\end{picture}

Here an unoriented edge can be oriented in an arbitrary way. As usual, a branch can even be replaced by a rooted tree with appropriate zero relations ( see \cite[ 10.7]{Buch} ). The same remark applies to all extended Dynkin-diagrams. All categories obtained from a Dynkin-diagram of type $T$ by orienting the edges and by replacing certain branches are then called of type $T$.

For  later use we collect  some simple facts in the following lemma.

\begin{lemma}\label{easy}
Let $P$ and $D$ be ray categories.
\begin{enumerate}
 \item Let $\sigma$ and $\rho$ be morphisms with in $P$ with common domain. If there are morphisms $\phi$ and $\psi$ such that $\sigma \phi =0 \neq \rho \phi$ and $\sigma \psi \neq 0= \rho \psi$, or  if there is a morphism $\eta$ such that $\rho\eta$  and $\sigma \eta$ are neighbors in a crown then the following diagram is cleaving:

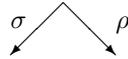
\begin{picture}(20,2)
 \put(8,2){\vector(-1,-1){1}}
\put(8,2){\vector(1,-1){1}}
\put(7,1.5){$\sigma$}
\put(9,1.5){$\rho$}
\end{picture}
\item The composition of cleaving functors is cleaving.
\item If $\tau$ is long in $P$ and $F:D \rightarrow P$ is cleaving with $F\mu \neq \tau$ for all $\mu$ in $D$, then the induced functor $D \rightarrow P/\tau$ is still cleaving.
\item  $F:D \rightarrow P$ is cleaving iff it satisfies the conditions: i) $F\mu=0$ iff $\mu=0$. ii) No irreducible morphism is mapped to an identity. iii) For any two irreducible morphisms $\alpha:x \rightarrow y$  and $\beta:x \rightarrow z$ in $D$ and each $\preceq$-maximal morphism  $\mu:x \rightarrow t$ with $\beta \preceq \mu$ that does not factor through $\alpha$, the image $F\mu$ does not factor through $F\alpha$. iv) The dual of iii).

\item The category $D$ given by the quiver in figure 3.1 without relations contains the crown $(\alpha \gamma,  \beta \gamma,  \beta \delta,  \alpha \delta
)$ of length $4$.  Similarly, the category $D$ given by the  quiver in figure 3.2 with zero-relations $\alpha_{2}\alpha_{1}$,$\beta_{2}\beta_{1}$ and $\gamma_{2}\gamma_{1}$ contains the crown $(\beta_{2}\alpha_{1},\beta_{2}\gamma_{1},\alpha_{2}\gamma_{1},\alpha_{2}\beta_{1},\gamma_{2}\beta_{1},\gamma_{2}\alpha_{1})$ of length $6$.\vspace{0.5cm}

\setlength{\unitlength}{0.8cm}
\begin{picture}(20,2)

\put(4,2){\vector(1,-1){1}}

\put(6,2){\vector(-1,-1){1}}

\put(5,1){\vector(-1,-1){1}}

\put(5,1){\vector(1,-1){1}}

\put(5.75,0.4){$\beta$}
\put(3.75,0.4){$\alpha$}
\put(3.75,1.4){$\gamma$}
\put(5.75,1.4){$\delta$}
\put(4,-1){figure 2.1}

\put(10,2){\vector(1,-1){1}}

\put(11,2){\vector(0,-1){1}}
\put(12,2){\vector(-1,-1){1}}

\put(11,1){\vector(-1,-1){1}}
\put(11,1){\vector(0,-1){1}}
\put(11,1){\vector(1,-1){1}}

\put(12,0.2){$\gamma_{2}$}
\put(9.4,0.2){$\alpha_{2}$}
\put(9.4,1.8){$\alpha_{1}$}
\put(11.1,1.8){$\beta_{1}$}
\put(11.1,0.2){$\beta_{2}$}
\put(12,1.8){$\gamma_{1}$}
\put(10,-1){figure 2.2}

\end{picture}\vspace{0.8cm}

\item Let $(v=\alpha_{s}\alpha_{s-1} \ldots \alpha_{1}, w= \beta_{t}\beta_{t-1},\ldots \beta_{1})$ be a contour in $P$. Then there is the cleaving diagram shown in figure 3.3.

\setlength{\unitlength}{1cm}
\begin{picture}(10,3)
\put(1,1){\vector(1,1){1}}
\put(2,2){\vector(1,0){1}}
\put(3,2){\vector(1,0){1}}
\put(4,2){......}
\put(5,2){\vector(1,0){1}}
\put(6,2){\vector(1,-1){1}}
\put(1,1){\vector(1,-1){1}}
\put(2,0){\vector(1,0){1}}
\put(3,0){\vector(1,0){1}}
\put(4,0){......}
\put(5,0){\vector(1,0){1}}
\put(6,0){\vector(1,1){1}}

 \put(1,1.5){$\alpha_{1}$}
\put(2.2,2.2){$\alpha_{2}$}
\put(5.2,2.2){$\alpha_{s-1}$}
\put(5.2,0.2){$\beta_{t-1}$}
\put(6.8,1.5){$\alpha_{s}$}
\put(6.8,0.5){$\beta_{t}$}
\put(2.2,0.2){$\beta_{2}$}

\put(1,0.5){$\beta_{1}$}
\put(2,0){\vector(1,0){1}}
\put(3,0){\vector(1,0){1}}
\put(4,0){......}
\put(5,0){\vector(1,0){1}}
\put(6,0){\vector(1,1){1}}\put(3,-1){figure 2.3}

\end{picture}

\end{enumerate}

\end{lemma}\vspace{1cm}
 
The easy proofs are left to the reader. The parts a) and its dual, c) and d) are used again and again in the technical sections to come.\vspace{0.8cm}

At the end of this paragraph I make some  comments on the fundamental article \cite{BGRS} on multiplicative bases. It consists of a local and a global part.

The local part deals only with small pieces of the given algebra $A$. Here one uses quite often the cleaving technique. Unfortunately, one has to deal with arbitrary $k$-categories instead of ray-categories which makes the verification of the cleaving conditions much more complicated. This is one reason why the local part is hard to read. However, due to \cite{Standard}, it suffices to deal only with $\vec{A}$ instead of $A$ which is much  easier. This considerable simplification mentioned already in \cite{BGRS} is explained with many details in chapter 13 of the book \cite{Buch}. In a forthcoming paper \cite{contours} we refine the structure and disjointness theorems for non-deep contours and  obtain shorter proofs especially for the so called diamonds.

After the rather technical local part there is the global topological part starting with section 8. This part is already very elegant and independent of lengthy case-by-case considerations. Nevertheless, in \cite{Fischbacher} it was simplified and generalized by Fischbacher whose main results are presented in  13.9 and 14.2 of the book \cite{Buch}. But observe that the crucial reduction lemma in \cite{Fischbacher} is based on a beautiful lemma about hooks in efficient tackles contained in \cite[ 8.4]{BGRS} resp.  \cite[lemma 25]{Bongo}).

\subsection{Long morphisms in crowns}

The next result is basic for the inductive proof of theorem 2.

\begin{proposition}\label{keyprop}Let $P$ be a minimal representation-infinite ray category containing at least one crown and one long morphism. Then there is a long morphism not occuring in a contour.
 
\end{proposition}

We will explain now the general strategy how to prove this and we give the details in the following sections.  To each crown $C$ we consider the pair of natural numbers $(n,t)$, where $2n$ is the length of the crown and $t$ is the sum $\sum _{i=1}^{n} (d(\rho_{i}) + d(\sigma_{i}))$ of the depths of all morphisms in $C$. The lexikographic order on the pairs $(n,t)$ induces a partial order on the set of crowns. We choose a minimal element $C$. Of course, each long morphism $\tau$ occurs in $C$, because otherwise $C$ induces a crown in $P/{\tau}$ contradicting the fact that $P$ is minimal representation-infinite. 

It is clear that we are in a self-dual situation: The minimal crown $C$ is also a minimal crown in the minimal representation-infinite ray-category $P^{op}$, that contains also a long morphism. Furthermore,  the proposition holds for $P$ iff it holds for $P^{op}$. So if we have proved that $C$ has a certain property, it has also the dual property.
   
Assume that the proposition is not true i.e. that all long morphisms belong to a contour. We will  derive in several steps the contradiction that $C$ as above does not exist. 

We can assume that $\sigma_{1}$ is long, and we choose a contour $$(v=\alpha_{s} \alpha_{s-1}\ldots \alpha_{1}, w=\beta_{t}\beta_{t-1}\ldots \beta_{1})$$ with $\vec{v}=\vec{w}=\sigma_{1}$. We say that $\rho_{1}$ factors through $v$ resp. $w$ if we have $\rho_{1}= \rho_{1}'\vec{\alpha_{1}}$ resp. $\rho_{1}=\overline{\rho_{1}} \vec{\beta_{1}}$.

 We show in 2.2.1 that $\rho_{1}$ factors through exactly one of the two paths $v$ or $w$.
Dually, one defines when $\rho_{n}$ factors through $v$ resp. $w$. Of course, by the above self-duality, $\rho_{n}$ also factors through exactly one of the two paths. So there are two cases possible for a contour $(v,w)$ belonging to the long morphism $\sigma_{1}$. Either both neighbours $\rho_{1}$ and $\rho_{n}$ factor through different paths or both factor through the same path. In the first case the chosen contour is called permeable, in the second reflecting.
Analogous definitions and statements hold for all long morphisms occurring in $C$ and for all choices of contours.

In the very technical section 2.2.2 we show that two long morphisms are not neighbors in $C$  and in 2.2.3 that $n=2$ is not possible.
In 2.2.4 we look at the long morphism $\sigma_{1}$ with the chosen contour $(v,w)$ and we assume that we have a non-trivial factorization  $v=v_{2}v_{1}$ such that $\rho_{1}=\rho_{1}'\vec{v}_{1}$. We prove that then $\sigma_{2}$ does not factor through $\rho_{1}'$. Finally in 2.2.5 we show that $C$ does not exist. 

For the proofs of the first and especially the last step one has to look at large parts of the category whereas the other proofs only require a careful analysis of some small parts. However this local part is more complicated than in the article on multiplicative bases because we have to consider also deep contours.

 We end this section with some easy, but useful observations.

\begin{lemma}\label{factorsatmost} We keep all the notations and assumptions made in this section.
\begin{enumerate}

\item For $x\epsilon P$, the ray category induced by the partially orderered set $S:=x/P$ is  zigzag-free. 

\item $\rho_{1}$ does not factor through $v$ and $w$.

\item Suppose $\rho_{1} \leq \sigma_{1}$. If $\vec{v}=\delta\gamma$ is a non-trivial factorization and $\rho_{1}=\rho_{1}'\gamma$, then $\delta$ is not irreducible. Furthermore, there is a non-trivial factorization $v=v_{3}v_{2}v_{1}$  such that $\rho_{1}=\rho_{1}'\vec{v}_{1}$ and such that the  diagram in figure 4.1 is cleaving. Here $w=w_{2}w_{1}$ is any non-trivial decomposition of $w$.

\setlength{\unitlength}{0.5cm}
\begin{picture}(9,4)
 \put(5,4){\vector(-1,-1){2}}
\put(3,2){\vector(1,-1){2}}
\put(5,4){\vector(1,-1){1}}
\put(6,3){\vector(1,-1){2}}
\put(6,3){\vector(0,-1){2}}
\put(6,1){\vector(-1,-1){1}}
\put(7.5,2){$\rho_{1}'$}
\put(3,3){$\vec{w}_{1}$}
\put(5,3){$\vec{v}_{1}$}
\put(3,0.7){$\vec{w}_{2}$}
\put(5.2,2){$\vec{v}_{2}$}
\put(5,0.8){$\vec{v}_{3}$}\put(4,-1){figure 4.1}

\end{picture}\vspace{0.5cm}

 \end{enumerate}
\end{lemma}

Proof: a) Suppose not. Since $S$ is finite, there is a crown $(\mu_{1},\nu_{1}, \ldots  , \nu_{m})$ in $S$. The morphism $\mu_{i}$ in $S$ gives us two morphisms $\phi_{i},\psi_{i}$ with domain $x$ and a morphism $\mu_{i}'$ satisfying $\psi_{i}=\mu_{i}'\phi_{i}$. Similarly, we obtain morphisms $\nu_{i}'$ with $\psi_{i+1}=\nu_{i}'\phi_{i}$, where $\psi_{m+1}=\psi_{1}$. Here no $\phi_{i}$ is the identity of $x$, because otherwise $\mu_{i}$ factors through $\nu_{i-1}$ ( $\nu_{0}:= \nu_{m}$ ).  It follows easily  that $(\mu_{1}',\nu_{1}', \ldots  , \nu_{m}')$ is a crown in $P$ that does not contain $\tau$, where $\tau$ is a non-zero morphism of maximal depth. This is impossible because $P$ is minimal representation infinite.

b) So assume $\rho_{1}$ factors through both, i.e. we have $\rho_{1}= \rho'\vec{\alpha_{1}}=\rho''\vec{\beta_{1}}$. Set $v'=\alpha_{s} \ldots \alpha_{2}$ and $w'=\beta_{t}\ldots \beta_{2}$. Then $C'=(\rho',\vec{v'},\vec{w'},\rho'')$ is a crown in $x/P$, where $x$ is the domain of $\rho_{1}$.

c) By assumption we have  $\rho_{1}\leq \sigma_{1}$, i.e. $\sigma_{1}= \psi \rho_{1}\phi$ for some morphisms $\phi$ and $\psi$.  The definition of a ray category implies that one of the following four equations relating $\psi$ and $\delta \rho_{1}'$ holds for some appropriate $\xi$.
If we have $\xi \psi=\delta \rho_{1}'$   we get $\xi \psi \vec{v}_{1}=\delta \rho_{1}'\gamma=\xi \sigma_{1}\neq 0$. Because $\sigma_{1}$ is long, $     \xi$ is an identity. But then we obtain $\sigma_{1}=\psi\vec{v}_{1}=\delta\rho_{1}'\gamma= \delta \rho_{1}$ contradicting the fact that $\sigma_{1}$ and $\rho_{1}$ are neighbors in a crown.

The case $ \psi\xi =\delta \rho_{1}'$ leads to $0\neq  \psi \vec{v}_{1}=\delta \rho_{1}'\phi\gamma = \psi\xi\phi\gamma$. Cancellation shows that $\xi$ is an identity and we are in the impossible first case.

The third possibility $\psi= \xi \delta \rho_{1}'$  contradicts again to the fact that $\sigma_{1}$ and  $\rho_{1}$ are neighbors. Thus we are in the case $\psi= \delta \rho_{1}'\xi$ where $\xi$ is not an identity because otherwise we are back in the third case. Since $\rho_{1}'$ is not an identity either we are done.

Now, choose $i$ maximal with the property that $\rho_{1}$ factors through $\vec{v}_{1}$ where $v_{1}=\alpha_{i} \ldots \alpha_{1}$. Then we have $i<s$ as shown before which gives the wanted  non-trivial factorization. The diagram is cleaving by the maximal choice of $i$.

\subsubsection{$\rho_{1}$  factors through $v$ or $w$}

\begin{lemma}\label{factors}We keep all notations and assumptions.
\begin{enumerate}
	\item $\rho_{1}$ factors through $v$ or $w$.
	\item At most one neighbor $\rho_{i}$ of the long morphism $\sigma_{1}$ satisfies $\rho_{i}\leq \sigma_{1}$.
         \item There are at least two different long morphisms in $C$.
         \item If $\sigma_{1}$ and $\rho_{1}$ are long and $\sigma_{2}\leq \sigma_{1}$, then $\sigma_{2}\leq \rho_{1}$.
\end{enumerate}

\end{lemma}

Proof: a) Suppose a) is not true. We start with the long morphism $       \sigma_{1}$ and we move ahead in the crown until we reach the next long morphism $\tau$ which might be $\sigma_{1}$ again. We consider first the case where $\tau$ occurs as some $\rho_{i}$. We choose a contour $(v'=\alpha_{s'}'\ldots \alpha_{1}',w'=\beta_{t'}',\ldots ,\beta_{1}')$ corresponding to the long morphism $\tau$.
 
If $\tau= \rho_{1}$ the first arrows in $v,w,v',w'$ are all different. This is impossible. Thus we have $i>1$. If $\sigma_{i}$ does not factor through $v'$ or $w'$, we consider the morphisms $\alpha_{1},\beta_{1},\rho_{1}, \ldots ,\sigma_{i},\alpha_{1}',\beta_{1}'$. They define a quiver of type $\tilde{D}_{m}$ as a  cleaving diagram in $P/\sigma_{1}$ as can be seen in figure 5.1.\vspace{1cm}

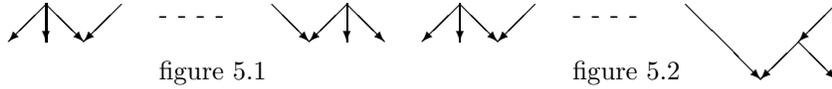
\setlength{\unitlength}{0.5cm}
\begin{picture}(22,2)

\put(1,2){\vector(-1,-1){1}}

\put(1,2){\vector(0,-1){1}}

\put(1,2){\vector(1,-1){1}}
\put(3,2){\vector(-1,-1){1}}\multiput(4,1.5)(0.5,0){4}{-}
\put(7,2){\vector(1,-1){1}}
\put(9,2){\vector(-1,-1){1}}
\put(9,2){\vector(1,-1){1}}
\put(9,2){\vector(0,-1){1}}
\put(4,0){figure 5.1}
\put(15,0){figure 5.2}
\put(12,2){\vector(-1,-1){1}}

\put(12,2){\vector(0,-1){1}}

\put(12,2){\vector(1,-1){1}}
\put(14,2){\vector(-1,-1){1}}\multiput(15,1.5)(0.5,0){4}{-}
\put(18,2){\vector(1,-1){2}}
\put(21,1){\vector(-1,-1){1}}
\put(22,2){\vector(-1,-1){1}}

\put(21,1){\vector(1,-1){1}}

\end{picture}\vspace{0.5cm}

So suppose that we have a non-trivial decomposition $v'=v_{2}'v_{1}'$ such that $\sigma_{i}=\sigma_{i}'\vec{v}_{1}'$. First we look at the subcase where $\rho_{i-1}$ does not factor through $\sigma_{i}'$ and we consider the cleaving diagram of type $\tilde{D}_{m}$ defined by the morphisms $\alpha_{1},\beta_{1},\rho_{1}, \ldots ,\sigma_{i}',\vec{v}_{1}',\vec{v}_{2}'$ and drawn in figure 5.2. If $\sigma_{1}\neq \tau$, this lies in $P/\sigma_{1}$ which is impossible. If $\sigma_{1}= \tau$ and if there is another long morphism $\phi$ we obtain the same cleaving diagram in $P/\phi$.
Finally if $\sigma_{1}= \tau$ is the only long morphism we have $\sigma_{i}\leq \rho_{i}$, whence there is by lemma \ref{factorsatmost} a non-trivial decomposition $v_{2}'=uu'$ with an arrow $u$. Replacing $\vec{v}_{2}'$ by $\vec{u}'$ we find again a $\tilde{D}_{m}$- quiver in $P/\tau$.

We are left with the case $\rho_{i-1}=\sigma_{1}'\rho_{i-1}'$. For $i-1=1$ one gets a $\tilde{D}_{5}$ -quiver in $P/\sigma_{1}$ defined by the morphisms $\alpha_{1},\beta_{1},\rho_{1}',\sigma_{1}',\vec{v}_{1}'$. For $i>2$ we go on with these factorizations  $\sigma_{i-1}= \sigma_{i-1}'\rho_{i-1}'$, $\rho_{i-2}= \rho_{i-2}'\sigma_{i-1}'$ and so on as long as possible. If all morphisms can be factorized we end up with a similar cleaving diagram as above with $\rho_{2}'$ instead of $\vec{v}_{1}'$.

If $\rho_{k}$ is the first morphism that does not factor  we obtain in $P/\sigma_{1}$ the  cleaving diagram of figure 5.3 defined by the morphisms $\alpha_{1},\beta_{1},\rho_{1}, \ldots ,\rho_{k},\sigma_{k+1}',\rho_{k+1}',\sigma_{k+2}'$, and if it is some $\sigma_{k}$ the similar cleaving diagram of figure 5.4.\vspace{1cm}

\setlength{\unitlength}{0.5cm}
\begin{picture}(20,2)\put (-1,0){\begin{picture}(13,2)

\put(1,2){\vector(-1,-1){1}}

\put(1,2){\vector(0,-1){1}}

\put(1,2){\vector(1,-1){1}}
\put(3,2){\vector(-1,-1){1}}\multiput(4,1.5)(0.5,0){4}{-}\put(5,-1){figure  5.3}
\put(8,2){\vector(1,-1){2}}
\put(11,1){\vector(-1,-1){1}}
\put(11,1){\vector(0,-1){1}}
\put(11,2){\vector(0,-1){1}}

\end{picture}}
\put(12,0){
\begin{picture}(13,2)

\put(1,2){\vector(-1,-1){1}}

\put(1,2){\vector(0,-1){1}}

\put(1,2){\vector(1,-1){1}}
\put(3,2){\vector(-1,-1){1}}\multiput(4,1.5)(0.5,0){4}{-}\put(5,-1){figure 5.4}
\put(8,2){\vector(-1,-1){2}}
\put(8,2){\vector(1,-1){1}}
\put(9,2){\vector(0,-1){1}}
\put(9,1){\vector(0,-1){1}}

\end{picture}}\end{picture}\vspace{1cm}

The case where $\tau$ is some $\sigma_{i}$ can be treated with the same arguments. 

b) Suppose that $\sigma_{1}$  goes from $x$ to $y$.  By  axiom e) of a ray category $P(x,y)$ is cyclic over $P(x,x)$ or over $P(y,y)$. If $P(x,y)$ is generated by $\chi$ over $P(x,x)$ we cannot have $\rho_{n}\leq \sigma_{1}$. For choosing a generator $\gamma$ of $P(x,x)$ we obtain from $\sigma_{1}=\psi \rho_{n}\phi$ the relations $\sigma_{1}=\chi \gamma^{s}$ and $\rho_{n}\phi = \chi \gamma^{t}$ with $s\geq t$, whence $\sigma_{1}=\rho_{n}\phi\gamma^{s-t}$. This contradicts the fact, that $\sigma_{1}$ and $\rho_{n}$ are neighbors in a crown. Dually, if $P(x,y)$ is cyclic over $P(y,y)$, we cannot have $\rho_{1}\leq \sigma_{1}$.

c) Assume that $\sigma_{1}$ is the only long morphism.
Then the neighbors $\rho_{1}$ and $\rho_{n}$ factor through $v$ or $w$ by part a), whence they are not irreducible. Therefore they are smaller than the only long morphism $\sigma_{1}$ which contradicts part b).

d) By assumption, we have $\sigma_{1}=\psi\sigma_{2}\phi\neq 0$. By axiom e) of a ray category we have $\sigma_{2}\phi \leq \rho_{1}$ and also $\sigma_{2}\leq \rho_{1}$.

\subsubsection{Long morphisms are not neighbors}

This lengthy section is only devoted to prove:

\begin{lemma}\label{notneighbor} Two long morphisms are not neighbors in $C$.
 
\end{lemma}

Proof: Suppose on the contrary that $\sigma_{1}$ and $\rho_{1}$ are long. For $\sigma_{1}$ we take the already chosen contour $(v,w)$ and for $\rho_{1}$ we choose an arbitrary contour $(v',w')$.  Let $x$ be the domain of $\sigma_{1}$. At most three arrows start at $x$. Thus  we can assume that $\alpha_{1}$ is the first arrow of $v$ and $v'$.  Inside the partially ordered set  $x/P$ we look at the set $S$ of all morphisms $\phi$ satisfying $\vec{\alpha}_{1}\preceq \phi$, $\phi \preceq \sigma_{1}$ and $\phi \preceq \rho_{1}$. This set contains a greatest element $\psi$ because otherwise $S$ contains a crown of length 2 contradicting part a) of lemma \ref{factorsatmost}. Now there are paths $u=\gamma_{r}\ldots \gamma_{1}$, $u'=\delta_{p}\delta_{p-1}\ldots \delta_{1}\gamma_{r'}\ldots \gamma_{1}$ and  $r'<r$ such that $\alpha_{1}=\gamma_{1}$, $\psi=\vec{\gamma}_{r'}\ldots \vec{\gamma}_{1}$, $\vec{u}=\sigma_{1}$  and $\vec{u'}=\rho_{1}$. Since $u$ and $v$ are interlaced by construction,  $(u,w)$ is a contour  with $\vec{u}=\sigma_{1}$. Similarly, $(u',w')$ is a contour with $\vec{u'}=\rho_{1}$. Now define $u_{2}=\gamma_{r}\ldots \gamma_{r'+1}$,$u'_{2}=\delta_{p}\delta_{p-1}\ldots \delta_{1}$. Then we obtain for any  non-trivial decompositions $w=w_{2}w_{1}$, $w'=w'_{2}w'_{1}$ the following cleaving diagram in $P$:

\setlength{\unitlength}{0.5cm}
\begin{picture}(18,6)

\put(10,6){\vector(1,-1){2}}
\put(7.5,5){$\vec{w}_{1}$}

\put(10,5){$\psi$}

\put(12,5){$\vec{w}_{1}'$}

\put(7,2.5){$\vec{w}_{2}$}

\put(12.5,2.5){$\vec{w}_{2}'$}

\put(9,2){$\vec{u}_{2}$}

\put(10.5,2){$\vec{u}_{2}'$}

\put(10,6){\vector(-1,-1){2}}

\put(10,6){\vector(0,-1){2}}
\put(8,4){\vector(0,-1){2}}
\put(10,4){\vector(-1,-1){2}}
\put(10,4){\vector(1,-1){2}}
\put(12,4){\vector(0,-1){2}}
\put(9,0){figure 6.0}

\end{picture}\vspace{0.5cm}

Namely,  $\sigma_{1}$ does not factor through $\vec{w}_{1}'$ by part b) of lemma \ref{factorsatmost}. By symmetry, $\rho_{1}$ does not factor through $\vec{w}_{1}$. Moreover $\vec{u}_{2}=\xi \vec{u}_{2}'$ implies $\sigma_{1}=\xi \rho_{1}$. This is impossible since $\sigma_{1}$ and $ \rho_{1}$ are neighbors in a crown. For the same reason, $\vec{u}'=\xi \vec{u}$ leads to a contradiction.

Observe that $\rho_{n}$ cannot factor through $\vec{w}_{2}$ and $\vec{u}_{2}$. The analogous statement holds for $\sigma_{2}$ and we have the tedious task to analyze the different possibilities. This is not difficult, but very lengthy. We give always a representation-infinite ray category that is cleaving by its number in figure 2.2 or by the type of an extended Dynkin diagram, but we do not check in detail all the conditions imposed on a cleaving functor. For instance, part a) of lemma \ref{easy} and its dual will be used very often without mentioning it explicitely.

First let $\sigma_{2}$ be long too. Then in one of the two paths of a contour to $\sigma_{2}$ there is an arrow $\eta$ with the same codomain as $\rho_{1}$ such that $\rho_{1}$ does not factor through $\vec{\eta}$.  If $\rho_{n}$ is also long one gets an  arrow $\theta$ having the corresponding properties with respect to $\sigma_{1}$. This gives in $P/\sigma_{1}$ the cleaving diagram from figure 6.1 involving from the left to the right the morphisms $\theta,\vec{w}_{2}, \vec{u}_{2},\vec{u}_{2}',\vec{w}_{2}',\eta$.

\setlength{\unitlength}{0.5cm}
\begin{picture}(20,7)

\put(6,3){\vector(-1,-1){1}}\put(5.5,1.4){$\vec{\eta}$}

\put(0,3){\vector(1,-1){1}}\put(0,1.4){$\vec{\theta}$}
\put(1,4){\vector(0,-1){2}}
\put(3,4){\vector(-1,-1){2}}
\put(3,4){\vector(1,-1){2}}
\put(5,4){\vector(0,-1){2}}\put(2,0.5){figure 6.1}

\put(16,3){\vector(-1,-1){1}}

\put(11,6){\vector(0,-1){2}}
\put(11,4){\vector(0,-1){2}}
\put(13,4){\vector(-1,-1){2}}
\put(13,4){\vector(1,-1){2}}
\put(15,4){\vector(0,-1){2}}\put(10,5){$\rho_{n}'$}

\put(11,4){\vector(-1,-1){2}}

\put(11,4){\vector(0,-1){2}}\put(12,0.5){figure 6.2}

\end{picture}

So $\rho_{n}$ is not long. If it does not factor through $\vec{w}_{2}$ or $\vec{v}_{2}$ we obtain a similar cleaving diagram as above with $\rho_{n}$ instead of $\vec{\theta}$. So suppose $\rho_{n}=\vec{w} _{2}\rho_{n}'$. If $\sigma_{n}=\sigma_{n}'\rho_{n}'$ we have in $P/\sigma_{1}$ the cleaving diagram drawn in figure 6.2 given by the morphisms $\rho_{n}',\sigma_{n}', \vec{w}_{2}, \vec{u}_{2},\vec{u}_{2}',\vec{w}_{2}',\eta$.

If $\sigma_{n}$ does not factor through $\rho_{n}'$ we get in $P/\sigma_{2}$ the cleaving diagram shown in figure 6.3. It involves the morphisms $\phi,\rho_{n}',\vec{w}_{2}, \vec{u}_{2},\vec{u}_{2}',\vec{w}_{2}',\vec{w}_{1}, \psi,\vec{w}_{1}',\eta$. Here $\phi$ is $\sigma_{n}$ for $\sigma_{n}\neq \sigma_{2}$ or else an irreducible morphism occuring in the paths chosen to a contour corresponding to the long morphism $\sigma_{2}$ that $\rho_{n}$ does not factor through.\vspace{0.5cm}

\setlength{\unitlength}{0.5cm}
\begin{picture}(20,7)

\put(6,6){\vector(1,-1){2}}

\put(6,6){\vector(-1,-1){2}}

\put(6,6){\vector(0,-1){2}}
\put(4,4){\vector(0,-1){2}}
\put(1.5,3){$\phi$}
\put(6,4){\vector(-1,-1){2}}
\put(6,4){\vector(1,-1){2}}
\put(8,4){\vector(0,-1){2}}
\put(9,3){\vector(-1,-1){1}}
\put(2,4){\vector(1,0){2}}
\put(2,4){\vector(0,-1){2}}
\put(4,0){figure 6.3}

\put(15,6){\vector(1,-1){2}}
\put(15,6){\vector(-1,-1){2}}
\put(15,6){\vector(0,-1){2}}
\put(13,4){\vector(0,-1){2}}
\put(15,4){\vector(-1,-1){2}}
\put(19,6){\vector(-1,-1){2}}
\put(14,0){figure 6.4}
\put(13,6){\vector(1,-1){2}}
\put(13,6){\vector(-1,-1){2}}
\end{picture}\vspace{0.5cm}

Next we consider the case where $\rho_{n}$ factors through $u$. We will always find an appropriate $\tilde{D}_{m}$- quiver that admits a cleaving functor into a proper quotient of $P$. If $\rho_{n}$ does not factor through $\vec{u}_{2}$ we have a nontrivial decomposition $u_{2}=ab$ with $\rho_{n}=\vec{b}\rho_{n}'$. By the definition of $\psi$, $\vec{u}_{2}'$ does not factor through $\vec{b}$ and the morphisms $\rho_{n}',\vec{b},\vec{a},\vec{u}_{2}',\vec{w}_{2}',\eta$ define our wanted cleaving functor into $P/\rho_{1}$. 

So we have $\rho_{n}=\vec{u}_{2}\rho_{n}'$. Suppose first, that $0\neq \vec{u}'_{2}\rho_{n}'$. Then $P$ contains the crown $(\rho_{n},\sigma_{1},\rho_{1},\vec{u}'_{2}\rho_{n}')$. Therefore we have $n=2$ and $\sigma_{2}=\vec{u}'_{2}\rho_{n}'$ by the minimal choice of $C$. By part d) of lemma \ref{factors} $\rho_{n}$ is comparable to $\sigma_{1}$ or to $\sigma_{2}$. From $\rho_{n}\leq \sigma_{1}$ we see using lemma \ref{factorsatmost} part c) that $\psi$ can be non-trivially factored as $\psi_{2}\psi_{1}$ such that the morphisms $\psi_{2},\vec{u}_{2},\vec{u}_{2}',\vec{w}_{2}',\eta$  give rise to a cleaving functor from a $\tilde{D}_{5}$-quiver to $P/\rho_{1}$.
Similarly, $\rho_{n}\leq \sigma_{2}$ implies that $u_{2}'$ admits a non-trivial factorization $u_{2}'=ab$. Then we find to $P/\rho_{1}$ the cleaving functor from a $\tilde{D}_{4}$- quiver given by the morphisms $\rho_{2}',\psi_{2},\vec{u}_{2},\vec{b}$.

 We are reduced to the case $0= \vec{u}'_{2}\rho_{n}'$. If $\sigma_{n}=\sigma_{n}'\rho_{n}'$ then we have in $P/\sigma_{2}$ the  $\tilde{D}_{5}$-diagram supported by the morphisms $\sigma_{n}',\vec{u}_{2},\vec{u}_{2}',\vec{w}_{2}',\eta$.

If this is not cleaving, we have $\xi\vec{u}'_{2}= \sigma_{n}'$ or else $\vec{u}'_{2}= \xi \sigma_{n}'$. The first case implies the contradiction $0=\xi\vec{u}'_{2}\rho_{n}'=\sigma_{n}'\rho_{n}'=\sigma_{n}\neq 0$. In the second case we can assume that $\xi$ is not an identity. From $0\neq \vec{u}_{2}'\psi=\xi \sigma_{n}'\psi$ we obtain $0\neq \sigma_{n}'\psi$. Then we have in $P$ the  crown $(\sigma_{n}'\psi,\sigma_{n},\rho_{n},\sigma_{1})$ which is strictly smaller than the given chain $C$ because the depth of $\rho_{1}$ is strictly greater than the depth of $\sigma_{n}'\psi$.

 So from now on $\sigma_{n}$ does not factor through $\rho_{n}'$.  If $\sigma_{2}$ factors through $w'$ we can choose the decomposition so that it factors already through $\vec{w}'_{2}$. Then we obtain in $P/\rho_{1}$ or for $\sigma_{n}=\rho_{1}$ in $P/\sigma_{2}$ the cleaving diagram from figure 6.4 that shows an algebra with number $11$ from the list. The occurring morphisms are $\sigma_{n},\rho_{n}',\vec{w}_{2},\vec{w}_{1},\psi,\vec{u}_{2},\vec{w}_{1}',\sigma_{2}'$.

So $\sigma_{2}$ factors through $u'$. If we even have $\sigma_{2}=\vec{u}'_{2}\sigma_{2}'$ we obtain in figure 6.5 the next cleaving diagram in $P/\sigma_{1}$ containing the morphisms $\rho_{n}',\vec{w}_{1},\psi,\vec{w}_{1}',\sigma_{2}'$.:\vspace{0.5cm}

\setlength{\unitlength}{0.5cm}
\begin{picture}(20,7)

\put(3,6){\vector(1,-1){2}}

\put(3,6){\vector(-1,-1){2}}

\put(3,6){\vector(0,-1){2}}
\put(1,6){\vector(1,-1){2}}
\put(5,6){\vector(-1,-1){2}}
\put(2,1){figure 6.5}

\put(3,6){\vector(-1,-1){2}}

\put(12,6){\vector(-1,-1){2}}
\put(12,6){\vector(0,-1){2}}
\put(10,4){\vector(0,-1){2}}
\put(12,4){\vector(-1,-1){2}}
\put(12,4){\vector(1,-1){1}}
\put(14,4){\vector(-1,-1){1}}

\put(12,6){\vector(1,-1){1}}\put(10.5,5.5){$\rho_{n}'$}\put(12.5,3.6){$\vec{b}$}

\put(10,6){\vector(1,-1){2}}

\put(10,6){\vector(-1,-1){1}}
\put(10,1){figure 6.6}\put(16,4){$\vec{b}\rho_{n}'=0$}

\end{picture}

 Here $\sigma_{2}'=\rho_{n}'\xi$ leads to $0\neq \vec{u}_{2}'\sigma_{2}'=\vec{u}_{2}'\rho_{n}'\xi=0$, whereas $\rho_{n}'=\sigma_{2}'\xi$ with a non-identity $\xi$ implies $0\neq \vec{u}_{2}\sigma_{2}'$. Then we find in $P$ a $\tilde{D}_{4}$-quiver as a cleaving diagram inducing the crown $(\sigma_{1},\vec{u}_{2}\sigma_{2}',\sigma_{2},\rho_{1})$ which is strictly smaller than $C$.

In the last case remaining with a long $\sigma_{2}$ there is a non-trivial decomposition  $u_{2}'=ab$ such that $\sigma_{2}=\vec{a}\sigma_{2}'$. Then we get $\vec{b}\rho_{n}'=0$ because otherwise we have an obvious cleaving diagram of trype $\tilde{D}_{4}$ in $P/\rho_{1}$. Now figure 6.6 shows a cleaving functor from the category  with number 20 from our list into $P/\rho_{1}$. The  occurring morphisms are $\sigma_{n},\rho_{n}',\vec{w}_{2},\vec{w}_{1},\psi,\vec{u}_{2},\vec{w}_{1}',\vec{b},\sigma_{2}'$. Recall here that $\sigma_{n}$ does not factor through $\rho_{n}'$.

 We have treated all cases where $\sigma_{2}$ is long. By duality, we are reduced to the situation that neither $\rho_{n}$ nor $\sigma_{2}$ are long. First, let $\rho_{n}$ factor through $w$. Choosing an appropriate decomposition for $w$ we have $\rho_{n}=\vec{w}_{2}\rho_{n}'$. Then $u'_{2}$ is an arrow, because for $u'_{2}=ab$ one has the  cleaving functor shown in figure 6.7 from the category with number 11 to $P/\rho_{1}$. The occurring morphisms are $\rho_{n}',\vec{w}_{2},\vec{w}_{1},\psi,\vec{u}_{2},\vec{w}_{1}',\vec{b}.$

If $\sigma_{2}$ factors through $w'$ we find a decomposition such that $\sigma_{2}=\vec{w}'_{2}\sigma_{2}'$. This gives  in $P$ the category  93 as a cleaving diagram drawn in figure 6.8 with morphisms $\rho_{n}',\vec{w}_{2},\vec{w}_{1},\psi,\vec{u}_{2},\vec{w}_{1}',\vec{u}_{2}',\vec{w}_{2}',\sigma_{2}'$.

\vspace{0.5cm}

\setlength{\unitlength}{0.5cm}
\begin{picture}(20,7)

\put(3,6){\vector(-1,-1){2}}

\put(3,6){\vector(0,-1){2}}
\put(1,4){\vector(0,-1){2}}
\put(3,4){\vector(-1,-1){2}}
\put(3,4){\vector(1,-1){1}}
\put(0,4){\vector(1,0){1}}

\put(3,6){\vector(1,-1){1}}
\put(1,0){figure 6.7}

\put(11,6){\vector(-1,-1){2}}

\put(11,6){\vector(0,-1){2}}
\put(9,4){\vector(0,-1){2}}
\put(11,4){\vector(-1,-1){2}}
\put(11,4){\vector(1,-1){2}}
\put(13,4){\vector(0,-1){2}}

\put(11,6){\vector(1,-1){2}}

\put(14,4){\vector(-1,0){1}}

\put(8,4){\vector(1,0){1}}
\put(10,0){figure 6.8}

\end{picture}\vspace{0.5cm}

If there is a third long morphism $\tau$ different from $\sigma_{1}$ and $\rho_{1}$ this diagram is in $P/\tau$. In the other case we have $\rho_{n}\leq \sigma_{1}$ and  $\sigma_{2}\leq \rho_{1}$ by lemma \ref{factors} part d). Thus part c) of lemma \ref{factorsatmost} implies that   $w_{1}=ab$ and $w'_{1}=a'b'$ with non-trivial decompositions and arrows $b,b'$. and we find in figure 6.9 a cleaving diagram from a  $\tilde{D}_{8}$-quiver  to $P/\sigma_{1}$. The occurring morphisms are $\rho_{n}',\vec{w}_{2},\vec{a},\vec{u}_{2},\vec{u}_{2}',\vec{w}_{2}',\vec{a}',\sigma_{2}'$. Observe that $\rho_{n}'$ resp. $\sigma_{n}'$ does not factor through $\vec{a}$ resp. $\vec{a}'$ by  lemma \ref{factorsatmost} again.\vspace{0.5cm}

\setlength{\unitlength}{0.5cm}
\begin{picture}(7,6)

\put(1,5){\vector(0,-1){1}}
\put(3,4){\vector(-1,-1){2}}
\put(3,4){\vector(1,-1){2}}
\put(5,5){\vector(0,-1){1}}
\put(5,4){\vector(0,-1){2}}
\put(1,4){\vector(0,-1){2}}

\put(6,4){\vector(-1,0){1}}

\put(0,4){\vector(1,0){1}}
\put(1.5,0){figure 6.9}

\put(11,6){\vector(-1,-1){2}}

\put(11,6){\vector(0,-1){2}}
\put(11,4){\vector(-1,-1){2}}
\put(9,4){\vector(0,-1){2}}
\put(13,6){\vector(-1,-1){2}}

\put(11,6){\vector(1,-1){2}}
\put(7,4){\vector(1,0){2}}

\put(9,0){figure 6.10}

\put(17,6){\vector(1,-1){2}}

\put(17,6){\vector(0,-1){2}}

\put(17,4){\vector(1,-1){2}}
\put(19,6){\vector(-1,-1){2}}
\put(19,4){\vector(0,-1){2}}
\put(19,6){\vector(1,-1){2}}

\put(17,6){\vector(-1,-1){2}}\put(20.2,5){$\phi$}

\put(14,4){\vector(1,0){1}}\put(16,0){figure 6.11}

\end{picture}\vspace{0.5cm}

If $\sigma_{2}$ factors through $u'$ it factors already through the irreducible morphism $\vec{u}_{2}'$, i.e. we have $\sigma_{2}=\vec{u}_{2}'\sigma_{2}'$. Then we have $\vec{u}_{2}\sigma_{2}'=0$ because otherwise $P/\rho_{1}$ contains a  category of type 11 as a cleaving diagram. This is  indicated in figure 6.10. The morphisms involved are $\rho_{n}',\vec{w}_{2},\vec{w}_{1},\psi,\vec{u}_{2},\vec{w}_{1}',\sigma_{2}'$.

If $\rho_{2}$ does not factor through $\sigma_{2}'$  we obtain the cleaving functor of figure 6.11 from a category of type 12 to  $P/\sigma_{1}$. Here $\phi$  is $\rho_{2}$ for $\rho_{2}\neq \sigma_{1}$ or else an irreducible morphism that belongs to the path $u$ or $w$ that $\sigma_{2}$ does not factor through. The second construction always works for $\rho_{2}=\sigma_{1}$. The invoved morphisms are $\rho_{n}',\vec{w}_{1},\psi,\vec{u}_{2}',\vec{w}_{1}',\vec{w}_{2}',\sigma_{2}',\phi$.

So we have $\rho_{2}= \rho_{2}'\sigma_{2}'$ and $\rho_{2}\neq \sigma_{1}$. If $\rho_{2}'\psi\neq 0$ we find a $\tilde{D}_{4}$ quiver consisting of the morphisms $\psi, \sigma_{2}',\rho_{2}',\vec{u}_{2}'$ in $P/\sigma_{1}$. Thus we have $\rho_{2}'\psi=0$ and we find in figure 6.12 an $\tilde{E}_{8}$-quiver or in figure 6.13  a $\tilde{D}_{6}$-quiver in $P/\sigma_{1}$ depending on the fact whether $\sigma_{n}$ factors through $\rho_{n}'$ or not. The involved morphisms are obvious.

\vspace{0.5cm}

\setlength{\unitlength}{1cm}
\begin{picture}(18,3)
\put(-1,0){
\begin{picture}(7,3)

\put(3,2){\vector(1,-1){1}}
\put(1,2){\circle*{0.1}}

\put(3,2){\vector(0,-1){1}}
\put(3,2.5){$\sigma_{2}'$}

\put(3,2){\vector(-1,-1){1}}
\put(4,2){\vector(0,-1){1}}
\put(3.2,1){$\rho_{2}'$}
\put(4,3){\vector(-1,-1){1}}
\put(2,2){\vector(0,-1){1}}

\put(1,2){\vector(1,0){1}}

\put(1,2){\vector(-1,0){1}}
\put(4.5,1){$\vec{u}_{2} \sigma_{2}'=0$}
\put(2,0.5){figure 6.12}

\end{picture}}
\put(6,0){
\begin{picture}(7,3)

\put(3,2){\vector(1,-1){1}}

\put(3,2){\vector(0,-1){1}}

\put(3,2){\vector(-1,-1){1}}
\put(2,2){\vector(0,-1){1}}
\put(2,2){\vector(-1,-1){1}}
\put(1,2){\vector(1,0){1}}\put(2,0.5){figure 6.13}

\end{picture}}
\end{picture}

By the symmetry of our situation we are left with the case where $\rho_{n}$ factors through $u$ and $\sigma_{2}$ through $u'$. One of the two paths $u_{2}$ or $u'_{2}$ is an arrow because otherwise we have an obvious $\tilde{D}_{5}$- quiver in $P/\sigma_{1}$. Say $u_{2}$ is an arrow. Then we have $\sigma_{2}=\vec{u}'_{2}\sigma_{2}'$.  If we have also $\rho_{n}=\vec{u}_{2}\rho_{n}'$ we look at the  diagram of figure 6.14 in $P/\sigma_{1}$. It contains the morphisms $\rho_{n}',\vec{w}_{1},\psi ,\vec{w}_{1}', \sigma_{2}'$. If this diagram is not cleaving we have up to symmetry $\sigma_{2}'= \rho_{n}'\xi$. This implies $\vec{u}'_{2}\rho_{n}' \neq 0$ and we obtain the crown $(\rho_{n},\sigma_{1},\rho_{1},\vec{u}'_{2}\rho_{n}')$. The minimal choice of $C$ implies $n=2$. From parts c) resp. d) of lemmata \ref{factorsatmost} resp. \ref{factors}  we get that $\rho_{n}\leq \sigma_{1}$ and  that $\psi$ has a non-trivial factorization $\psi= \psi_{2}\psi_{1}$ with irreducible $\psi_{1}$. Then figure 6.15 shows  a $\tilde{D}_{4}$-quiver in $P/\sigma_{1}$ consisting of the morphisms $\rho_{n}',\vec{u}_{2},\psi_{2},\vec{u}_{2}'$.
\vspace{0.5cm}

\setlength{\unitlength}{0.5cm}
\begin{picture}(20,4)

\put(3,4){\vector(1,-1){2}}

\put(3,4){\vector(0,-1){2}}

\put(3,4){\vector(-1,-1){2}}
\put(1,4){\vector(1,-1){2}}
\put(5,4){\vector(-1,-1){2}}\put(1.5,0){figure 6.14}

\put(10,4){\vector(0,-1){1}}

\put(9,4){\vector(1,-1){1}}
\put(10,3){\vector(-1,-1){1}}
\put(10,3){\vector(1,-1){1}}
\put(8.5,0){figure 6.15}

\put(17,4){\vector(1,-1){2}}

\put(17,4){\vector(0,-1){2}}

\put(17,4){\vector(-1,-1){1}}
\put(16,3){\vector(-1,-1){1}}
\put(15,3){\vector(1,0){1}}
\put(15.5,0){figure 6.18}

\end{picture}\vspace{0.5cm}

Finally, $\rho_{n}$ does not factor through $\vec{u}_{2}$. Then $u_{2}$ has a non-trivial decomposition $u_{2}=ba$ and $\rho_{n}=\vec{b}\rho_{n}'$. One has $\vec{a}\sigma_{2}'=0$ because otherwise a $\tilde{D}_{4} $-quiver with morphisms $\psi,\sigma_{2}',\vec{a},\vec{u}_{2}'$ is cleaving in $P/\sigma_{1}$. If $\sigma_{n}$ does not factor through $\rho_{n}'$ we have in $P/\sigma_{1}$ an $\tilde{E}_{6}$-quiver as shown in figure 6.16.\vspace{0.5cm}

\setlength{\unitlength}{0.5cm}
\begin{picture}(20,4)
\put(2,-1.5){figure 6.16}

\put(5,2){\vector(-1,-1){2}}
\put(7,4){\vector(-1,-1){2}}
\put(7,4){\vector(1,-1){2}}
\put(3,4){\vector(0,-1){4}}
\put(1,2){\vector(0,-1){2}}
\put(1,2){\vector(1,0){4}}

\put(16,4){\vector(1,-1){2}}

\put(16,4){\vector(0,-1){2}}

\put(16,2){\vector(1,-1){2}}
\put(18,2){\vector(0,-1){2}}
\put(16,4){\vector(-1,-1){2}}
\put(16,2){\vector(-1,-1){2}}
\put(12,2){\vector(1,-1){2}}
\put(18,4){\vector(-1,-1){2}}
\put(18,4){\vector(1,-1){2}}
\put(17.5,3){$\sigma_{2}'$}
\put(15.5,0.8){$\vec{a}$}
\put(19,0){$\vec{a}\sigma_{2}'=0$}
\put(15.5,-1.5){figure 6.17}

\end{picture}\vspace{1cm}

If $\rho_{2}$ does not factor through $\sigma_{2}'$ we find in $P/ \sigma_{1}$ a category of type 20 as a cleaving diagram. This is illustrated in figure 6.17. So suppose we have $\rho_{2}=\rho_{2}'\sigma_{2}'$. Then  $P/\sigma_{1}$ contains the $\tilde{D}_{5}$-quiver shown in figure 6.18. as a  diagram. The morphisms involved are $\rho_{n}',\vec{b}, \vec{a}, \rho_{2}',\vec{u}_{2}'$. If this diagram is not cleaving, we have $\vec{b}\vec{a}=\xi\rho_{2}'$ with a non-identity $\xi$. But then $P/\sigma_{1}$ contains as  cleaving diagram the $\tilde{D}_{4}$-quiver with morphisms $\psi,\sigma_{2}',\rho_{2}',\vec{u}_{2}'$.

\subsubsection{The length of $C$ is at least $6$}

\begin{lemma}\label{nnottwo} The length of $C$ is not $4$.
\end{lemma}

Proof:  Suppose $n=2$. By part c) of lemma \ref{factors}, there are at least two different long morphisms which cannot be neighbors by the last section. So up to duality we have by part b) of lemma \ref{factors} that $\sigma_{1}$ and $\sigma_{2}$ are long and that  $\rho_{i}\leq \sigma_{i}$ holds for $i=1,2$. So by part c) of lemma \ref{factorsatmost} we obtain from $\sigma_{1}$ and $\rho_{1}$ a cleaving diagram as in figure 2.3. Choosing a contour $(v',w')$ for $\sigma_{2}$ we obtain a similar diagram from $\sigma_{2}$ and $\rho_{2}$ thereby assuming that $\rho_{2}$ factors through $v'$. Then there is an arrow $\eta$ in $v'$ or in $w'$ that $\rho_{1}$ does not factor through. So we obtain in $P/\sigma_{2}$ the cleaving diagram shown in figure 7.1. If one of the $v_{i}$'s is not an arrow or if the length of $w$ is $\geq 5$, then we obtain obviously an $\tilde{E}_{6}$- or an $\tilde{E}_{8}$-quiver as a cleaving diagram in $P/\sigma_{2}$. This is indicated in figure 7.2. \vspace{0.5cm}

\setlength{\unitlength}{0.5cm}
\begin{picture}(20,6)

 \put(2,4){\vector(-1,-1){2}}
\put(6,1){$\vec{\eta}$}
\put(6,2){\vector(-1,-1){1}}
\put(0,2){\vector(1,-1){2}}
\put(2,4){\vector(1,-1){1}}
\put(3,3){\vector(1,-1){2}}
\put(3,3){\vector(0,-1){2}}
\put(3,1){\vector(-1,-1){1}}
\put(4.5,2){$\rho_{1}'$}
\put(0,3){$\vec{w}_{1}$}
\put(3,3.5){$\vec{v}_{1}$}
\put(0,0.7){$\vec{w}_{2}$}
\put(2,2){$\vec{v}_{2}$}
\put(3,0.3){$\vec{v}_{3}$}
\put(1,-1){figure 7.1}

 \put(8,4){\vector(1,-1){1}}
\put(9,3){\vector(0,-1){1}}
\put(9,2){\vector(0,-1){1}}
\put(9,1){\vector(-1,-1){1}}
\put(9,2){\vector(1,-1){1}}
\put(11,2){\vector(-1,-1){1}}

 \put(12,4){\vector(1,-1){1}}
\put(12,4){\vector(-1,-1){1}}
\put(13,3){\vector(0,-1){1}}
\put(13,3){\vector(1,-1){1}}
\put(13,2){\vector(0,-1){1}}
\put(15,3){\vector(-1,-1){1}}

 \put(17,4){\vector(1,-1){1}}
\put(17,4){\vector(-1,-1){1}}
\put(18,3){\vector(0,-1){1}}
\put(18,3){\vector(1,-1){1}}
\put(16,3){\vector(0,-1){1}}
\put(16,2){\vector(0,-1){1}}
\put(16,1){\vector(0,-1){1}}
\put(20,3){\vector(-1,-1){1}}
\put(12,-1){figure 7.2}

\end{picture}\vspace{1cm}

Thus we can assume that all $v_{i}$ and all $v_{i}'$ are arrows and also that $(rad P)^{5}=0$. Otherwise there is a path $u$ of length $\geq 5$ such that $\vec{u}$ is a long morphism say $\sigma_{1}$ and we can replace the contour $(v,w)$ by $(v,u)$ or by $(u,w)$ which is impossible.

First we consider the case where the contour $(v',w')$ chosen to $\sigma_{2}$ is permeable. For $w_{1}$ we take an arrow and we factorize $\rho_{1}=\rho_{1}' \vec{v}_{1}$. If $\rho_{1}'$  factors through $w'$ we find the ray category 44 from our list as a cleaving diagram in $P/\sigma_{2}$ as drawn in figure 7.3.  The morphisms involved are $\vec{w}_{1},\vec{w}_{2},\vec{v}_{1},\vec{v}_{2},\vec{v}_{3},\rho_{1}',\vec{v}_{3}',\vec{v}_{2}',\rho_{2}'$. Here $\rho_{1}$ does not factor through $\vec{v}_{3}'$, because it factors already through $w'$. Also $\vec{v}_{3}'\vec{v}_{2}'$ does not factor through $\rho_{1}'$, because $v'$ and $w'$ are not interlaced.
\vspace{0.5cm}

\setlength{\unitlength}{0.5cm}
\begin{picture}(20,5)

\put(1,4){\vector(1,-1){1}}

\put(1,4){\vector(-1,-2){1}}

\put(0,2){\vector(1,-2){1}}
\put(2,3){\vector(0,-1){2}}
\put(2,1){\vector(-1,-1){1}}
\put(2,3){\vector(2,-3){2}}
\put(3,2){$\rho_{1}'$}
\put(5,1){\vector(-1,-1){1}}
\put(5,3){\vector(0,-1){2}}
\put(5,3){\vector(1,-1){1}}
\put(2,-1){figure 7.3}

\put(12,4){\vector(1,-1){1}}

\put(12,4){\vector(-1,-2){1}}

\put(15,4){\vector(1,-1){1}}

\put(12,4){\vector(2,-1){2}}\put(13,4){$\gamma$}

\put(15,4){\vector(-1,-1){1}}
\put(16,3){\vector(0,-1){2}}
\put(16,3){\vector(1,-1){1}}
\put(13,-1){figure 7.4}

\end{picture}\vspace{1cm}

If $\rho_{1}'$ does not factor through $w'$ we have another factorization $\rho_{1}= \vec{w}_{2}'\gamma$ giving rise in $P/\sigma_{1}$ to the  cleaving diagram shown in figure 7.4. It involves the morphisms $\vec{w}_{1},\vec{v}_{1},\gamma, \vec{w}_{1}',\vec{v}_{1}',\vec{v}_{2}',\rho_{2}'$.
Here $\gamma$ does not factor through $\vec{w}_{1}$, because $\rho_{1}$ does not factor through $w$, and it does not factor through $\vec{v}_{1}$, because $\rho_{1}'$ does not factor through $w'$.

By symmetry, we can assume now that the contours chosen to $\sigma_{1}$ and $\sigma_{2}$ are both reflecting. Then $\rho_{1}$ factors already  through $\vec{v}_{3}'\vec{v}_{2}'$, because for $\rho_{1}=\vec{v}_{3}'\phi$ we find a $\tilde{D}_{5}$-quiver as a cleaving diagram in $P/\sigma_{1}$ consisting of the morphisms $\phi,\vec{v}_{1}',\vec{v}_{2}',\vec{v}_{3}',\rho_{2}'$. By symmetry we are in the situation of figure 7.5.\vspace{0.5cm}

\begin{picture}(20,6)

 \put(7,4){\vector(-1,-1){2}}

\put(5,2){\vector(1,-1){2}}
\put(7,4){\vector(1,-1){1}}
\put(8,3){\vector(0,-1){2}}
\put(8,4){$\tilde{\rho}_{1}$}
\put(5,3){$\vec{w}_{1}$}
\put(7,3){$\vec{v}_{1}$}
\put(7.2,2){$\vec{v}_{2}$}
\put(8,1){\vector(-1,-1){1}}

\put(12,4){\vector(1,-1){2}}
\put(7,4){\vector(4,-1){4}}
\put(12,4){\vector(-4,-1){4}}

\put(14,2){\vector(-1,-1){2}}
\put(12,4){\vector(-1,-1){1}}
\put(11,3){\vector(0,-1){2}}
\put(10.5,4){$\tilde{\rho}_{2}$}
\put(13.5,3){$\vec{w}_{1}'$}
\put(11.5,3){$\vec{v}_{1}'$}
\put(11.2,2){$\vec{v}_{2}'$}
\put(11,1){\vector(1,-1){1}}

\put(8.5,-1){figure 7.5}

\end{picture}\vspace{1cm}

Here we have chosen $w_{1}$ and $w_{1}'$ as arrows. We have $\rho_{1}=\vec{v}_{3}'\vec{v}_{2}'\tilde{\rho}_{1}$ and $\rho_{2}=\vec{v}_{3}\vec{v}_{2}\tilde{\rho}_{2}$. Because $(rad P)^{5}=0$ the $\tilde{\rho}_{i}$'s are irreducible and $( \tilde{\rho}_{1},\vec{v}_{1}',\tilde{\rho}_{2},\vec{v}_{1})$ cannot be a crown in $P/\sigma_{1}$. This implies $\vec{v}_{1}=\tilde{\rho}_{1}$ and $\vec{v}_{1}'=\tilde{\rho}_{2}$.

 For $v_{1}=v_{1}'$ we get $v_{2}\neq v_{2}'$ because otherwise $\rho_{1}=\rho_{1}'\vec{v}_{2}\vec{v}_{1}$. Since $\rho_{1}\leq \sigma_{1}$ and $v_{3}$ is an arrow, this contradicts part c) of lemma \ref{factorsatmost}. Furthermore $w_{1}=w_{1}'$ implies $\rho_{1}=\vec{v}_{3}'\vec{v}_{2}'\vec{v}_{1}=\vec{w}_{2}' \vec{w}_{1}'=\vec{w}_{2}' \vec{w}_{1}$, i.e. $\rho_{1}$ factors through $w$. Thus we find a $\tilde{D}_{5}$-quiver as a cleaving diagram in $P/\sigma_{1}$ that involves the morphisms $\vec{w}_{1},\vec{v}_{1},\vec{w}_{1}', \vec{v}_{2},\vec{v}_{2}'$.

For $v_{1}\neq v_{1}'$ we get as above that $v_{2}\neq v_{2}'$ because $v_{3}$ is irreducible. But now one has a $\tilde{D}_{4}$- quiver in $P/\sigma_{1}$ consisting of the morphisms $\vec{v}_{1},\vec{v}_{2},\vec{v}_{1}',\vec{v}_{2}'$.

\subsubsection{The factorization of neighbors of long morphisms stops after one step}

We keep all the notations and use all the reductions already obtained. In particular, long morphisms are not neighbors and $n>2$. So there is no cleaving diagram as in figure 2.1 in $P$.
\begin{lemma}\label{facstop}
Let $(v,w)$ be a contour with $\vec{v}=\sigma_{1}$. Let $v=v_{2}v_{1}$ be any non-trivial decomposition with $\rho_{1}=\rho_{1}'\vec{v}_{1}$. Then $\sigma_{2}$ does not factor through $\rho_{1}'$.
 
\end{lemma}

Proof: Suppose on the contrary that $\sigma_{2}=\rho_{1}'\sigma_{2}'$. Then we have $\vec{v}_{2}\sigma_{2}'=0$ because otherwise $\vec{v}_{1},\vec{v}_{2},\sigma_{2}',\rho_{1}'$ define a cleaving diagram as in figure 2.1. First, we treat the case where the contour is permeable. If $\rho_{2}$ does not factor through $\sigma_{2}'$ we obtain  an obvious $\tilde{E}_{7}$-quiver drawn in figure 8.1 as a cleaving  diagram in $P/\tau$. Here we take $\tau=\sigma_{1}$ for $\sigma_{1}\neq \rho_{2}$ or else a long morphism different from $\sigma_{1}$.\vspace{0.5cm} \setlength{\unitlength}{1cm}
\begin{picture}(20,2)
 \put(1,2){\vector(1,-1){1}}

\put(3,2){\vector(1,-1){1}}
\put(3,2){\vector(-1,-1){1}}
\put(4,1){\vector(1,-1){1}}
\put(5,2){\vector(-1,-1){1}}
\put(5,2){\vector(1,-1){1}}
\put(7,2){\vector(-1,-1){1}}
\put(1,1.4){$\rho_{n}'$}
\put(1.9,1.4){$\vec{w}_{1}$}
\put(3,1.4){$\vec{v}_{1}$}
\put(4,0.4){$\rho_{1}'$}
\put(4.65,1.4){$\sigma_{2}'$}
\put(5.65,1.4){$\rho_{2}$}
\put(6.65,1.4){$\sigma_{3}$}
\put(4.5,-1){figure 8.1}

\end{picture}\vspace{1.5cm}

For $\rho_{2}=\rho_{2}' \sigma_{2}'$ we get $\rho_{2}'\vec{v}_{1}=0$ because figure 2.1 is not cleaving in $P$. 
If $\sigma_{n}=\sigma_{n}' \rho_{n}'$, figure 8.2 shows a cleaving diagram in $P/\sigma_{1}$.\vspace{0.5cm}

\setlength{\unitlength}{1cm}
\begin{picture}(20,2)

\put(4,1){\vector(-1,-1){1}}

\put(3,2){\vector(1,-1){1}}
\put(4,1){\vector(1,-1){1}}

\put(6,1){\vector(-1,-1){1}}
\put(6,1){\vector(1,-1){1}}
\put(6,1){\vector(3,-1){3}}
\put(4.5,-0.5){figure 8.2}

\put(6,0.4){$\rho_{1}'$}
\put(4,0.4){$\vec{w}_{2}$}
\put(4.9,0.4){$\vec{v}_{2}$}
\put(8.65,0.4){$\rho_{2}'$}
\put(3,1.4){$\rho_{n}'$}
\put(3,0.5){$\sigma_{n}'$}
\end{picture}\vspace{1cm}

If $\sigma_{3}$ does not factor through $\rho_{2}'$ we look at the cleaving diagram of type $\tilde{E}_{8}$ in $P$ that is given in figure 8.3. Here $\phi$
is $\sigma_{n}$ if this is not long or else an irreducible morphism that $\rho_{n}'$ does not factor through. For $\sigma_{1}\neq \rho_{2}$ this diagram is already in $P/\sigma_{1}$, while for $\sigma_{1}=\rho_{2}$ there is another long morphism $\tau$ and the diagram is in $P/\tau$.

If $\sigma_{3}=\rho_{2}'\sigma_{3}'$, it follows $\rho_{1}'\sigma_{3}' =0$ since figure 2.1 is not cleaving. Then the diagram of figure 8.4 is cleaving in $P$ where $\phi$ is defined as in the case before. Again this diagram is in $P/\sigma_{1}$ for $\sigma_{1}\neq \rho_{2}$ or else in $P/\tau$ for a long morphism $\tau\neq \sigma_{1}$.\vspace{0.5cm}

\setlength{\unitlength}{0.7cm}
\begin{picture}(20,2)

\put(1,1){\vector(1,-1){1}}

\put(3,1){\vector(1,-1){1}}

\put(4,2){\vector(-1,-1){1}}
\put(3,1){\vector(2,-1){2}}
\put(6,1){\vector(-1,-1){1}}
\put(3,1){\vector(-1,-1){1}}
\put(0,2){\vector(-1,-1){1}}
\put(0,2){\vector(1,-1){1}}

\put(2.7,0.4){$\rho_{1}'$}
\put(3.65,1.4){$\sigma_{2}'$}
\put(5.65,0.4){$\sigma_{3}$}
\put(4.65,0.4){$\rho_{2}'$}
\put(-1,0){$\vec{v}_{2}\sigma_{2}'=0$}
\put(2,-1){figure 8.3}
\put(-1,1.4){$\phi$}
\put(1,1.4){$\rho_{n}'$}
\put(1.8,0.4){$\vec{v}_{2}$}

\put(11,2){\vector(1,-1){1}}

\put(12,1){\vector(1,-1){1}}

\put(13,2){\vector(-1,-1){1}}
\put(12,1){\vector(2,-1){2}}
\put(14,2){\vector(-2,-1){2}}
\put(11,2){\vector(-1,-1){1}}
\put(9,2){\vector(-1,-1){1}}
\put(9,2){\vector(1,-1){1}}
\put(10.5,-1){figure 8.4}
\put(11.8,0.4){$\rho_{1}'$}
\put(11.85,1.4){$\sigma_{2}'$}
\put(10.85,1.4){$\vec{v}_{1}$}
\put(7,0){$\rho_{2}'\vec{v}_{1}=\rho_{1}'\sigma_{3}'=0$}
\put(13.65,0.4){$\rho_{2}'$}
\put(13.65,1.4){$\sigma_{3}'$}
\put(8,1.4){$\phi$}
\put(8.8,1.4){$\rho_{n}'$}

\end{picture}\vspace{1cm}

Now we consider the case where the contour $(v,w)$ is reflecting. Let $\tau\neq \sigma_{1}$ be another long morphism. We claim that $\rho_{n}=\vec{v}_{2}\rho_{n}'$. Suppose this is false. Then we get a proper decomposition $v_{2}=ba$ with $\rho_{n}=\vec{b}\rho_{n}'$. If $\rho_{1}'$ does not factor through $\vec{a}$, we have a $\tilde{D}_{5}$-quiver in $P/\tau$ involving the morphisms $\vec{v}_{1},\vec{a}, \vec{b},\rho_{n}',\rho_{1}'$. Thus we have $\rho_{1}'=\rho_{1}''\vec{a}$ with $\rho_{1}''\rho_{n}'=0$, because otherwise we obtain a $\tilde{D}_{4}$-quiver defined by the morphisms $\vec{a}, \vec{b},\rho_{n}',\rho_{1}''.$ If $\rho_{2}$ is long or if it does not factor through $\sigma_{2}'$, then $P/\sigma_{1}$ contains a cleaving diagram of type $\tilde{E}_{6}$. It involves the morphisms $\vec{w}_{1},\vec{v}_{1},\vec{a},\rho_{1}'',\sigma_{2}',\phi$, where $\phi$ is an arrow that $\sigma_{2}'$ does not factor through if $\rho_{2}$ is long or else $\phi=\rho_{2}$. Therefore, $\rho_{2}=\rho_{2}'\sigma_{2}'$. In case $\rho_{2}'=\rho_{2}''\vec{a}$ one has $\rho_{2}''\vec{a}\vec{v}_{1}=0$ since otherwise $P$ contains figure 2.1 as a cleaving diagram with the morphisms $\vec{a}\vec{v}_{1},\vec{a} \sigma_{2}',\rho_{1}'',\rho_{2}''$. But now the $\tilde{D}_{4}$-quiver defined by $\vec{a}, \vec{b},\rho_{1}'',\rho_{2}''$ is cleaving in $P$. Thus $\rho_{2}'$ does not factor through $\vec{a}$ and we also get $\rho_{2}'\vec{v}_{1}=0$ by the 'figure 2.1' argument. So we finally obtain in  $P/\sigma_{1}$ the $\tilde{E}_{7}$-quiver involving the morphisms $\vec{w}_{1},\vec{v}_{1},\vec{a}, \rho_{n}',\rho_{1}'',\sigma_{2}',\rho_{2}'$ with the relations $\rho_{1}''\rho_{n}'=\rho_{2}'\vec{v}_{1}=0$. We have shown our claim $\rho_{n}=\vec{v}_{2}\rho_{n}'$. Looking at the morphisms $\vec{v}_{1}, \rho_{n}',\rho_{1}',\vec{v}_{2}$ we infer $\rho_{1}'\rho_{n}'=0$. Then the diagram of figure 8.5 is cleaving in $P$.

\setlength{\unitlength}{0,7cm}
\begin{picture}(20,3)
 \put(3,2){\vector(0,-1){1}}
\put(1,1){\vector(1,-1){1}}
\put(2,2){\vector(1,-1){1}}
\put(2,2){\vector(-1,-1){1}}
\put(3,1){\vector(1,-1){1}}
\put(3,2){\vector(0,-1){1}}
\put(1,1){\vector(1,-1){1}}
\put(2,2){\vector(1,-1){1}}
\put(2,2){\vector(-1,-1){1}}
\put(3,1){\vector(1,-1){1}}
\put(3,1){\vector(-1,-1){1}}
\put(0.8,0.4){$\vec{w}_{2}$}
\put(1.8,0.4){$\vec{v}_{2}$}

\put(4,2){\vector(-1,-1){1}}
\put(2.5,2){$\rho_{n}'$}
\put(0.7,1.4){$\vec{w}_{1}$}
\put(1.9,1.4){$\vec{v}_{1}$}
\put(3.7,0.4){$\rho_{1}'$}
\put(3.7,1.4){$\sigma_{2}'$}\put(5,0){$\rho_{1}'\rho_{n}'=\vec{v}_{2}\sigma_{2}'=0$}
\put(2,-1){figure 8.5}

\put(11,2){\vector(0,-1){1}}
\put(10,2){\vector(1,-1){1}}
\put(10,2){\vector(-1,-1){1}}
\put(11,1){\vector(1,-1){1}}
\put(10,-1){figure 8.6}

\put(12,2){\vector(-1,-1){1}}
\put(12,2){\vector(1,-1){1}}
\put(10.4,2){$\rho_{n}'$}
\put(9.8,1.4){$\vec{v}_{1}$}
\put(10.7,0.4){$\rho_{1}'$}

\put(11.65,1.4){$\sigma_{2}'$}
\put(13,1.4){$\phi$}
\put(13.5,0){$\rho_{1}'\rho_{n}'=0$}

\end{picture}\vspace{1cm}

Define $\phi=\rho_{2}$ if $\rho_{2}$ is not long and does not factor through $\sigma_{2}'$ or take for $\phi$ an irreducible morphism that $\sigma_{2}$ does not factor through if $\rho_{2}$ is long. Then the $\tilde{E}_{6}$-diagram of figure 8.6 is cleaving in $P/\sigma_{1}$. Thus we have that $\rho_{2}$ is not long and $\rho_{2}=\rho_{2}'\sigma_{2}'$. This implies $\rho_{2}'\vec{v}_{1}=0$.

Suppose now that $\sigma_{n}$ does not factor through $\rho_{n}'$. Then $P$ admits an $\tilde{E}_{8}$-quiver as the cleaving diagram shown in figure 8.7. This lies already in $P/\sigma_{1}$ if $\sigma_{1}\neq \rho_{n-1}$. For $\sigma_{1}=\rho_{n-1}$ there is another long morphism $\tau$. For $\tau\neq \sigma_{2}$ the diagram lies in $P/\tau$. For $\tau=\sigma_{2}$ one gets an $\tilde{E}_{6}$-quiver drawn in figure 8.8 as a cleaving diagram in $P/\sigma_{2}$. Here $\phi$ is an irreducible morphism that $\sigma_{2}'$ does not factor through.

\vspace{1cm}

\setlength{\unitlength}{0.6cm}
\begin{picture}(20,3)

\put(2,2){\vector(-1,-1){1}}

\put(2,2){\vector(1,-1){1}}
\put(3,2){\vector(0,-1){1}}
\put(3,1){\vector(1,-1){1}}

\put(4,2){\vector(-1,-1){1}}
\put(4,2){\vector(1,-1){1}}
\put(6,2){\vector(-1,-1){1}}
\put(6,2){\vector(1,-1){1}}
\put(3,0){$\rho_{1}'$}
\put(2.4,2){$\sigma_{2}'$}
\put(2,1){$\vec{v}_{1}$}
\put(3.25,1){$\rho_{n}'$}
\put(6,0){$\rho_{1}'\rho_{n}'=0$}
\put(4.25,2){$\sigma_{n}$}
\put(5.25,1){$\rho_{n-1}$}
\put(6.25,2){$\sigma_{n-1}$}
\put(3.5,-1){figure 8.7}

\put(12,2){\vector(-1,-1){1}}

\put(12,2){\vector(1,-1){1}}
\put(13,3){\vector(0,-1){2}}
\put(13,3){\vector(1,0){1}}

\put(14,2){\vector(-1,-1){1}}
\put(14,2){\vector(1,-1){1}}
\put(13.4,2.5){$\phi$}
\put(12.4,2){$\sigma_{2}'$}
\put(12,1){$\vec{v}_{1}$}
\put(13.25,1){$\rho_{n}'$}
\put(14.25,2){$\sigma_{n}$}

\put(12,-1){figure 8.8}

\end{picture}\vspace{1.5cm}
Thus we can assume that $\sigma_{n}=\sigma_{n}'\rho_{n}'$ which implies $\sigma_{n}'\vec{v}_{1}=0$. Then the three morphisms $\vec{v}_{2},\rho_{1}',\rho_{2}'$ induce a cleaving diagram in $P$. If we add $\sigma_{n}'$ this cannot stay a cleaving diagram because none of the involved morphisms  is long. Thus we get $\xi \rho_{2}'=\sigma_{n}'$ or  $ \rho_{2}'=\xi \sigma_{n}'$. In both cases we obtain an $\tilde{E}_{6}$-quiver as a cleaving diagram in $P$. In the first case it involves the morphisms $\vec{v}_{1},\rho_{n}',\sigma_{2}',\vec{v}_{2},\rho_{1}',\rho_{2}'$ and in the second case $\vec{v}_{1},\rho_{n}',\sigma_{2}',\vec{v}_{2},\sigma_{n}',\rho_{1}'$. These induce by part e) of lemma \ref{easy} the two crowns $(\sigma_{1},\rho_{1},\sigma_{2},\rho_{2},\rho_{2}'\rho_{n}',\rho_{n})$ and $(\sigma_{1},\rho_{1},\sigma_{2},\sigma_{n}'\sigma_{2}',\sigma_{n},\rho_{n})$.
By minimality we get $n=3$ and $\xi=id$ in both cases. Since $\rho_{2}$ is not long, only the $\sigma_{i}$'s can be long. If all of them are long, we can always choose an irreducible morphism where $\vec{v}_{1}$ resp. $\rho_{3}'$ resp. $\sigma_{2}'$ does not factor through. This gives an $\tilde{E}_{6}$-quiver in a proper quotient $P/\tau$ as indicated in figure 8.9. Thus, using duality, we can assume that $\sigma_{1}$ and $\sigma_{2}$ are the only long morphisms. Then we get $\rho_{1}\leq \sigma_{1}$ or $\rho_{1}\leq \sigma_{2}.$ In the first case $v_{2}=ba$ by part c) of lemma \ref{factorsatmost}.  We find a $\tilde{D}_{4}$-quiver with morphisms $\vec{v}_{1},\rho_{3}',\sigma_{2}',\vec{a}$ or an $\tilde{E}_{6}$-quiver in $P/\sigma_{1}$ as shown in figure 8.10. The  case $\rho_{1}\leq \sigma_{2}$ is dual because there is also a reflecting contour to $\sigma_{2}$.

\vspace{0.5cm}

\setlength{\unitlength}{0.6cm}
\begin{picture}(20,5)

\put(2,2){\vector(-1,-1){1}}

\put(2,2){\vector(1,-1){1}}
\put(3,4){\vector(0,-1){3}}
\put(3,4){\vector(1,-1){1}}

\put(4,2){\vector(-1,-1){1}}
\put(4,2){\vector(1,-1){1}}
\put(2.3,2.5){$\rho_{3}'$}
\put(2.2,1){$\vec{v}_{1}$}
\put(3.6,1){$\sigma_{2}'$}
\put(4,-1){figure 8.9}

\put(12,2){\vector(-1,-1){1}}

\put(12,2){\vector(1,-1){1}}
\put(13,2){\vector(0,-1){1}}
\put(13,1){\vector(0,-1){1}}

\put(14,2){\vector(-1,-1){1}}
\put(13,1){\vector(1,-1){1}}
\put(14.2,0){$\rho_{1}'$}
\put(10.2,1){$\vec{w}_{1}$}
\put(12.2,1){$\vec{v}_{1}$}
\put(14.2,2){$\sigma_{2}'$}
\put(12.2,2){$\rho_{3}'$}
\put(12.5,0){$\vec{a}$}\put(15,1){$\rho_{1}'\rho_{3}'=\vec{a}\sigma_{2}'=0$}

\put(14,-1){figure 8.10}
\end{picture}\vspace{1cm}

\subsubsection{$C$ does not exist}

Assume that the contour $(v,w)$ chosen for the long morphism $\sigma_{1}$ in $C$ is permeable. Since the factorization of the neighbors stops after one step, we obtain up to permutation of $v$ and $w$ the following cleaving diagram drawn in figure 9.1 in $P$.\vspace{0.5cm}

\setlength{\unitlength}{1cm}\setlength{\unitlength}{1cm}
\begin{picture}(9,2)
 \put(2,2){\vector(1,-1){1}}
\put(2,2){\vector(-1,-1){1}}
\put(7,1){\vector(-1,-1){1}}
\put(3,1){\vector(1,-1){1}}
\put(4,2){\vector(1,-1){1}}
\put(4,2){\vector(-1,-1){1}}
\put(5,1){\vector(1,-1){1}}
\put(5,1){\vector(-1,-1){1}}

\put(2,1.4){$\rho_{n}'$}
\put(1,1.4){$\sigma_{n}$}
\put(7,0.4){$\sigma_{2}$}
\put(2.9,1.4){$\vec{w}_{1}$}
\put(4,1.4){$\vec{v}_{1}$}
\put(5,0.4){$\rho_{1}'$}
\put(2.9,0.4){$\vec{w}_{2}$}
\put(4,0.4){$\vec{v}_{2}$}\put(4,-1){figure 9.1}

\end{picture}\vspace{1.5cm}

Replacing $\sigma_{1}$ by $v_{2}$ and $\rho_{1}$ by $\rho_{1}'$ we obtain a smaller crown than $C$ because the sum of the depths of the involved morphisms has strictly decreased. 

We are left with the case where all the contours chosen for the long morphisms in the crown $C$ are reflecting. We start with $\sigma_{1}$ and go on until we reach the next long morphism $\tau$. Because long morphisms are not neighbors the first possible case is $\tau= \sigma_{2}$. We can assume that $\rho_{1}$ factors through $v$ and $v'$ where $(v',w')$ is a contour  with $\vec{v'}=\sigma_{2}$. Thus we  have non-trivial factorizations $v=v_{2}v_{1}$ and $\rho_{1}=\rho_{1}'\vec{v}_{1}$. If $\rho_{1}'$ factors through $v'$, we have $\rho_{1}'= \vec{v}_{2}'\rho_{1}''$ for some non-trivial factorization  $v'=v_{2}'v_{1}'$. Here $\rho_{1}''$ is not an identity because otherwise  $\sigma_{2}$ factors through $\rho_{1}'$ contradicting the last section. Thus we get the following representation infinite cleaving diagram shown in figure 9.2 in $P$:\vspace{0.5cm}

\setlength{\unitlength}{0.7cm}
\begin{picture}(14,6)
 \put(1,4){\vector(-1,-1){1}}
\put(1,4){\vector(1,-1){1}}
\put(0,3){\vector(1,-2){1}}
\put(2,3){\vector(-1,-2){1}}
\put(2,3){\vector(1,-1){1}}
\put(3,2){\vector(1,-1){1}}
\put(5,2){\vector(-1,-1){1}}
\put(4,4){\vector(-1,-2){1}}
\put(4,4){\vector(1,-2){1}}\put(2,-1){figure 9.2}

\put(10,4){\vector(-1,-1){1}}
\put(10,4){\vector(1,-1){1}}
\put(11,4){\vector(0,-1){1}}
\put(12,4){\vector(-1,-1){1}}
\put(12,5){\vector(0,-1){1}}
\put(12,5){\vector(-1,-1){1}}
\put(11,4){\vector(-1,-2){1}}
\put(9,2){\vector(1,0){1}}\put(9,-1){figure 9.3}

\end{picture}\vspace{1cm}

Because $P$ is minimal representation-infinite, $\sigma_{1}$ and $\sigma_{2}$ are the only long morphisms. But then $\rho_{1}\leq \sigma_{1}$ up to duality. By part c) of lemma \ref{factorsatmost}  $v_{2}=ba$ with an arrow $b$, and there is or a $\tilde{D}_{5}$-quiver as a cleaving diagram in $P/\sigma_{1}$. The morphisms involved are $\vec{v}_{1},\vec{a},\rho_{1}'',\vec{v}_{1}',\vec{v}_{2}'$. in the first case or in $P/\sigma_{2}$ in the second.

If $\rho_{1}'$ does not factor through $v'$, we obtain  in $P/\sigma_{1}$ the cleaving diagram shown in figure 9.3. It is a category with number 12 from our list.

Therefore we can assume that the minimal distance of two long morphisms in the crown is 3 at least. Next we consider the case where $\tau$ is some $\rho_{i}$. We have $1<i<n$. Using the last subsection we find the following cleaving diagram in $P$:\vspace{1.5cm}

\setlength{\unitlength}{1cm}
\begin{picture}(14,5)
\put(1,6){\vector(1,-3){1}}
\put(1,6){\vector(-1,-3){1}}
\put(0,3){\vector(1,-3){1}}
\put(2,3){\vector(-1,-3){1}}
\put(3,5){\vector(-1,-2){1}}
\put(2,3){\vector(1,-2){1}}
\put(3,5){\vector(1,-1){1}}
\put(5,5){\vector(-1,-1){1}}
 \put(4,2){\vector(-1,-1){1}}
\put(4,2){\vector(1,-1){1}}
\put(9,5){\vector(-1,-1){1}}
\put(9,5){\vector(1,-2){1}}
\put(11,6){\vector(1,-3){1}}
\put(11,6){\vector(-1,-3){1}}
\put(12,3){\vector(-1,-3){1}}
\put(10,3){\vector(1,-3){1}}

\put(10,3){\vector(-1,-2){1}}
\put(8,2){\vector(1,-1){1}}
\put(3.3,1){$\sigma_{2}$}
\put(1,3){$\sigma_{1}$}
\put(11,3){$\rho_{i}$}
\put(0.2,5){$\vec{w}_{1}$}
\put(1.5,5){$\vec{v}_{1}$}
\put(11.5,5){$\vec{w}_{1}'$}
\put(0.2,1){$\vec{w}_{2}$}
\put(8,1){$\rho_{i-1}$}
\put(11.5,1){$\vec{w}_{2}'$}
\put(2.5,5){$\rho_{n}'$}
\put(1.5,1){$\vec{v}_{2}$}
\put(6,4.5){.............}
\put(6,1.5){.............}
\put(9.2,5){$\sigma_{i+1}'$}
\put(10,5){$\vec{v}_{1}'$}
\put(8,5){$\rho_{i+1}$}
\put(10,1){$\vec{v}_{2}'$}
\put(9.2,1){$\sigma_{i}'$}
\put(4.5,1){$\rho_{2}$}
\put(3.3,5){$\sigma_{n}$}
\put(2.5,1){$\rho_{1}'$}
\put(4.2,0){$\rho_{1}'\rho_{n}'=0=\sigma_{i}'\sigma_{i+1}'$}
\put(4.2,5){$\rho_{n-1}$}\put(5,-1){figure 9.4}

\end{picture}\vspace{1.5cm}

Here we have denoted the contour chosen to the long morphism $\rho_{i}$ by $(v',w')$ and we have assumed that both neighbors factor through $v'$. The only zero-relations are $0=\rho_{1}'\rho_{n}'$ and $0=\sigma_{i}'\sigma_{i+1}'$. The two outer 'diamonds' connected by the lower finite zigzag form a representation-infinite cleaving diagram $D$ in $P$ that contains by construction only the two long morphisms $\sigma_{1}$ and $\rho_{i}$. Since $P$ is minimal representation-infinite, we conclude that these two morphisms are different and are the only long morphisms. If $v_{1}=ab$ is a non-trivial factorization, one gets an $\tilde{E}_{8}$-quiver as a cleaving diagram in $P/\sigma_{1}$ that involves the morphisms $\vec{a},\vec{w}_{2},\vec{v}_{2},\rho_{1}',\sigma_{2},\rho_{2},\sigma_{3},\rho_{3}$. The same argument shows that $v_{2}$,$v_{1}'$ and $v_{2}'$ are also arrows. If we had $\rho_{1} \leq \sigma_{1}$ we would obtain from part c) of lemma \ref{factorsatmost}, that $\vec{v}_2$ is not irreducible. Because $\rho_{1}=\rho_{1}' \vec{v}_{1}$ is neither irreducible nor long it is  comparable to a long morphism. In our situation we obtain $\rho_{i}=\delta \rho_{1}'\vec{v}_{1}\gamma$ for some appropriate morphisms $\gamma$ and $\delta$.

We claim that $(\vec{v}_{1}\gamma,\rho_{n}',\sigma_{n},\ldots  ,\rho_{i+1},\sigma_{i+1}',\vec{v}_{1}')$ is a crown. Since it is strictly smaller than $C$ this is impossible. Only four  factorisations between neighbors in the chain remain to be excluded. The factorisation $\vec{v}_{1}\gamma = \rho_{n}'\eta$  implies $0\neq \delta \rho_{1}' \vec{v}_{1}\gamma = \delta \rho_{1}'\rho_{n}'\eta =0$. From $\vec{v}_{1}\gamma\eta = \rho_{n}'$ we obtain 
$ \rho_{n}=\vec{v}_{2}\rho_{n}'=\vec{v}_{2}\vec{v}_{1}\gamma\eta =\sigma_{1}\gamma\eta $ which is ompossible for the neighbors $\sigma_{1}$ and $\rho_{n}$. Next assume $\eta \vec{v}_{1}'= \vec{v}_{1}\gamma$. We get $0\neq \delta \rho_{1}'\vec{v}_{1}\gamma = \delta \rho_{1}' \eta \vec{v}_{1}'= \rho_{i}= \vec{v}_{2}'\vec{v}_{1}'$, whence $\delta \rho_{1}' \eta = \vec{v}_{2}' $. Since all $v_{i}$ and $v_{i}'$ are arrows we see that $\eta$,$\delta$  and $\gamma$ are identities and that $\vec{v}_{2}'=\rho_{1}'$, $ \vec{v}_{1}'= \vec{v}_{1}$. This implies the contradiction $\rho_{1}=\rho_{i}$. Finally, look at $\vec{v}_{1}'= \eta \vec{v}_{1}\gamma$. This implies $v_{1}=v_{1}'$. In particular, the domains of the two long morphisms in $P$ coincide. 

Now we prove the same for the case where the next long morphism after $\sigma_{1}$ is some $\sigma_{i}$.  We can assume in addition $2<i<n$. Then we obtain a similar cleaving diagram as before:\vspace{1.5cm}

\setlength{\unitlength}{1cm}
\begin{picture}(14,5)
\put(1,6){\vector(1,-3){1}}
\put(1,6){\vector(-1,-3){1}}
\put(0,3){\vector(1,-3){1}}
\put(2,3){\vector(-1,-3){1}}
\put(3,5){\vector(-1,-2){1}}
\put(2,3){\vector(1,-2){1}}
\put(3,5){\vector(1,-1){1}}
\put(5,5){\vector(-1,-1){1}}
 \put(4,2){\vector(-1,-1){1}}
\put(4,2){\vector(1,-1){1}}
\put(9,5){\vector(-1,-4){1}}
\put(9,5){\vector(1,-2){1}}
\put(11,6){\vector(1,-3){1}}
\put(11,6){\vector(-1,-3){1}}
\put(12,3){\vector(-1,-3){1}}
\put(10,3){\vector(1,-3){1}}

\put(10,3){\vector(-1,-2){1}}
\put(8,5){\vector(1,-4){1}}
\put(3.3,1){$\sigma_{2}$}
\put(1,3){$\sigma_{1}$}
\put(11,3){$\sigma_{i}$}
\put(0.2,5){$\vec{w}_{1}$}
\put(1.5,5){$\vec{v}_{1}$}
\put(11.5,5){$\vec{w}_{1}'$}
\put(0.2,1){$\vec{w}_{2}$}
\put(7.2,1){$\sigma_{i-1}$}
\put(11.5,1){$\vec{w}_{2}'$}
\put(2.5,5){$\rho_{n}'$}
\put(1.5,1){$\vec{v}_{2}$}
\put(6,4.5){.............}
\put(6,1.5){.............}
\put(9.2,5){$\rho_{i-1}'$}
\put(10,5){$\vec{v}_{1}'$}
\put(7.2,5){$\sigma_{i+1}$}
\put(10,1){$\vec{v}_{2}'$}
\put(9.2,1){$\rho_{i}'$}
\put(4.5,1){$\rho_{2}$}
\put(3.3,5){$\sigma_{n}$}
\put(2.5,1){$\rho_{1}'$}
\put(4.2,5){$\rho_{n-1}$}\put(4.2,0){$\rho_{1}'\rho_{n}'=0=\rho_{i}'\rho_{i-1}'$}
\put(5,-1){figure 9.5}

\end{picture}\vspace{1.5cm}

Observe that now the two zigzags connecting the outer diamonds cross each other. The only zero-relations are $0=\rho_{1}'\rho_{n}'$ and $0=\rho_{i}'\rho_{i-1}'$. Argueing as  before one reduces to the situation where $\sigma_{1}$ and $\sigma_{i}$ are the only long morphisms, $\vec{v}_{1},\vec{v}_{2},\vec{v}_{1}'$ and $\vec{v}_{2}'$ are irreducible, $\rho_{1}$ is smaller than $\sigma_{i}$ and $\rho_{i}$ smaller than $\sigma_{1}$. So we have for some appropriate morphisms the equations $\delta \rho_{1}'\vec{v}_{1}\gamma = \sigma_{i}$ and $\delta '\rho_{i}'\vec{v}_{1}'\gamma' = \sigma_{1}$. 

We claim that $(\vec{v}_{1}\gamma, \rho_{n}',\sigma_{n}, \ldots \sigma_{i+1},\rho_{i})$ is a crown in $P$. Again only four factorizations between neighbors are not yet excluded by obvious reasons. $\rho_{n}'=\vec{v}_{1}\gamma \eta$ implies $\rho_{n}=\vec{v}_{2}\rho_{n}'=\sigma_{1} \gamma \eta$ which is impossible for neighbors in a crown.  Applying $\delta\rho_{1}'$ to the  equation $\rho_{n}'\eta =\vec{v}_{1}\gamma$ leads to $0=\delta\rho_{1}'\rho_{n}'\eta=\delta\rho_{1}'\vec{v}_{1}\gamma \neq 0$.
The third factorization to be excluded is $\eta \rho_{i}'\vec{v}_{1}'=\vec{v}_{1}\gamma.$ This gives $\delta \rho_{1}'\vec{v}_{1}\gamma = \sigma_{i}=\delta \rho_{1}'\eta \rho_{i}'\vec{v}_{1}'=\delta \rho_{1}'\eta \rho_{i}$ contradicting the fact that $\rho_{i}$ and $\sigma_{i}$ are neighbors in a crown. Finally suppose $ \rho_{i}'\vec{v}_{1}' = \eta  \vec{v}_{1}\gamma$. We get $\sigma_{1}=\delta '\rho_{i}'\vec{v}_{1}'\gamma'=\delta '\eta \vec{v}_{1}\gamma\gamma'$ so that $\delta' \eta \neq 0$. So we have $\delta' \eta = \vec{v}_{2}\xi$ or $\delta' \eta = \xi \vec{v}_{2}.$ In the second case $\xi$ is an identity because we have $0\neq \delta '\eta \vec{v}_{1}=\xi \vec{v}_{2} \vec{v}_{1}=\xi \sigma_{1}$. So we only have to consider the first case. From $0 \neq \vec{v}_{2} \vec{v}_{1}= \sigma_{1}= \delta '\eta \vec{v}_{1}\gamma\gamma'= \vec{v}_{2}\xi \vec{v}_{1}\gamma\gamma'$ we see that $\gamma$, $\gamma'$ and $\xi$ are all identities. In particular, the domains of the two long morphisms coincide again.

We have shown that $P$ contains only two long morphisms that have the same domain $x$. By the self-duality of our situation they also have the same codomain $y$. This contradicts the fact that in any ray-category $P(x,y)$ is linearly ordered.

\section{The proofs of the main results}
\subsection{The proof of theorem 2}

We need the following easy lemma.
\begin{lemma}
 Let $P$ be a ray category containig a long morphism $\mu$ that does not occur in a contour. Then we have:
\begin{enumerate}
 \item $P$ and $P/\mu$ have the same quivers and the same contours.
\item The quivers of the universal covers and the fundamental groups coincide.
\item For all abelian groups $Z$ the cohomology group ( \cite[section 8]{BGRS} ) $H^{2}(P,Z)$ embeds into $H^{2}(P/\mu ,Z)$.
\end{enumerate}

\end{lemma}

Proof: a) The quivers coincide because $\mu$ is long, and the contours, because $\mu$ does not belong to a contour.

b) The quivers of the universal covers and the fundamental groups are defined by the homotopy relation on the universal covers of the common quiver of $P$ and $P/\mu $. ( see \cite[ section 14.1]{Buch} or \cite[ section 10]{BGRS}). Since the definition of homotopic walks depends only on the contours, part b) follows from part a).

c) As shown in  \cite[section 8.2]{BGRS}, $H^{2}(P,Z)$ is isomorphic to the quotient of the space $C(P,Z)$ of $Z$-valued contour functions by the space $E(P,Z)$ of exact contour functions. The point is that here one has to take the 'old' definition of contour as given in \cite[section 2.7]{BGRS}, where our contours are called essential contours. Now in the old language part a) says that $P$ and $P/\mu$ have the same quivers and the same essential contours, but $P$ has in general more contours than $P/\mu$. Restriction is a homomorphism from $C(P,Z)$ to $C(P/\mu,Z)$ which induces an isomorphism on the spaces of exact contour functions because the quivers coincide. Since a contour function is uniquely determined by its values on the essential contours which are the same for $P$ and $P/\mu$, the restriction induces the wanted embedding.\vspace{0.5cm}

Now, let $A$ be a distributive minimal representation-infinite algebra. Then there is the associated ray category $P$ (\cite[ section 1.7]{BGRS} or  \cite[section 13.4]{Buch}). By theorem 13.17 in \cite{Buch}, that is based on \cite{BGRS} and \cite{Standard}, $P$ is also minimal representation-infinite. To prove the theorem we distinguish three cases.

If there is no long morphism in $P$, the ray category and the algebra are given by zero-relations of length 2, whence there are no contours.  Then theorem 2 is a well-known fact from elementary algebraic topology.

If there is no crown in $P$, then $P$ is zigzag-free in the terminology of \cite{Buch} and we are done by  \cite[theorems 14.2,13.17a]{Buch}. 

Finally, in the last case  there is a crown and a long morphism. By proposition 2, there is a long morphism $\mu$ as in the lemma above. Since $P/\mu$ contains no crown, the universal cover of $P/\mu$ is interval-finite and the fundamental group is free. By the lemma, the same holds for $P$. Theorem 13.17 a) of \cite{Buch} says that $A$ is isomorphic to $k^{f}(P)$ for some cohomology class in $H^{2}(P,k^{\ast})$. But this cohomology group vanishes for $P/\mu$ and by the lemma also for $P$. Therefore, $A$ is isomorphic to $k(P)$.
\subsection{The proof of theorem 1}

We have to show that a basic algebra  of infinite representation type has indecomposable representations in all dimensions.

 If $A$ is not distributive, we obtain this from proposition \ref{nondisinf}. If $A$ is distributive, we can assume that it is minimal representation-infinite and therefore by theorem 2 isomorphic to the linearization $k(P)$ of its ray-category. Moreover, the fundamental group is free. Thus the dimension  preserving push-down functor associated to the universal cover preserves indecomposability by ( \cite[section14.4]{Buch} or \cite{Galoiscover} ). So it is enough to find indecomposables of all dimensions for the universal cover $\tilde{P}$.

Section 2.3 in \cite{Bretscher} shows that the first homology group $H_{1}(k(\tilde{P}))$ of the Schurian category $k(\tilde{P})$ vanishes. Because the partially ordered sets $x/P$ and $P/y$ are zigzag-free by part a) of lemma 3, it follows from \cite[section 2.3]{Criterion} that each finite convex subcategory $B$ of $\tilde{P}$ satisfies also $H_{1}B=0$. By the separation-criterion \cite{Larrion} of Bautista-Larri\'{o}n as slightly generalized in \cite[section2.5]{Criterion}
any such $B$ has a preprojective component in its Auslander-Reiten quiver. 

We have to distinguish two cases.

 If there is such a subcategory $B$ that is not representation-finite, we apply some results of Ringel  in \cite[section 4.3]{Tame}. Namely there is a quotient $C$ of $B$ that is tame concealed. In particular, $C$ has indecomposable representations in all dimensions. Note that in this case by \cite[section 6]{Degenerations} each non-simple indecomposable is again an extension of an indecomposable and a simple as in the representation-finite case.

If all these finite subcategories are representation-finite there can be no common bound for the dimensions of all the indecomposable representations of $\tilde{P}$. For then the push-down functor would produce a bounded component of the Auslander-Reiten quiver of $P$ by \cite[ theorem 14.4]{Buch}. But then the category $P$, which is connected being minimal representation-infinite, would be representation-finite by a basic result in \cite[chapter 6]{Auslander}. So there are indecomposables of arbitrarily large dimensions, whence indecomposables in all dimensions by the known representation-finite case.

Observe that the proof of theorem 1 uses no classification lists at all. In contrast, all proofs of the second Brauer-Thrall conjecture via coverings depend on the list of the large faithful simply connected algebras in  \cite{Treue} ( but not on the lists in \cite{Listebild} or \cite{liste} as indicated in \cite[chapter 6 ]{Auslander} ).

\subsection{The proof of corollary 1}
Let $U$ be an indecomposable in $C$ of length $n$ and height $h$. Let $C(h)$ be the full subcategory of $C$ consisting of objects of height at most $h$, that have only the composition factors of $U$ as simple subquotients and that are of finite length. Then $C(h)$ is an abelian subcategory containing the indecomposable $U$. It is well-known that $C(h)$ is a module category (\cite{Unz}, \cite[section 8]{Indecom}) over some algebra which is finite-dimensional if all extension-groups between simples are finite-dimensional. For the convenience of the reader we give some details. 

If one of the extension groups $Ext(S,T)$ between simples in $C(h)$ is not finite-dimensional, one constructs easily local modules  of arbitrary length $n>2$ having top $S$. Indeed one takes $n-1$ linearly independent elements $E_{1},E_{2}, \ldots E_{n-1}$ in $Ext(S,T)$ and looks at the exact sequence $E:0 \rightarrow T^{n-1} \rightarrow X \rightarrow S \rightarrow 0$ such that the push-out under the projection $\pi_{i}:T^{n-1} \rightarrow T$ is $E_{i}$. Then $X$ is the wanted local module.

If all extension groups between the simples are of finite dimension one constructs finitely many projective indecomposables $P_{i}(h)$ whose direct sum is a progenerator $P(h)$ of finite length inside $C(h)$. We proceed  by induction on $h$. For $h=1$ we set $P_{i}(1)=S_{i}$ for a representative system $S_{1},S_{2},\ldots ,S_{r}$ of the  composition factors of $U$. In the inductive step we set $ dim_{k} Ext(P_{i}(h-1),S_{j}) = n_{ij}$. Note that by the half-exactness of $Ext$ these extension-groups are finite-dimensional. We define $P_{i}(h)$ as the ( uniquely determined ) middle term of the universal extension
$$0 \longrightarrow \bigoplus_{j=1}^{r} S_{j}^{n_{ij}} \longrightarrow P_{i}(h) \longrightarrow P_{i}(h-1) \longrightarrow 0.$$ We leave it as an exercise to show  that $P_{i}(h)$ is the projective cover of $S_{i}$ in $C(h)$. The functor $Hom(P(h), \,)$ identifies $C(h)$ with the finite dimensional right modules  over the finite dimensional endomorphism algebra of $P(h)$. Thus the corollary follows from theorem 1.
\subsection{The proof of corollary 2}

The implication from a) to b) is trivial. Reversely, theorem 1 implies that all indecomposables have length at most $n$. By Roiters theorem in \cite{BTI} the algebra is representation-finite. 

The bound $2 \cdot dim A + 1000$ is given in \cite[ section 5]{Treue}. It can be refined to the maximum of $2\cdot dim A$ and $30$.

\end{document}